\documentclass[journal]{IEEEtran}
\usepackage{cite}
\usepackage[hidelinks]{hyperref}
\usepackage{amsmath,amssymb,amsfonts}
\usepackage{tabularx}
\usepackage{comment}
\usepackage{booktabs}
\usepackage{threeparttable}
\usepackage{multirow}
\usepackage{siunitx}
\usepackage{relsize}
\usepackage{textcomp}
\usepackage{xcolor}
\usepackage{amsfonts,amsmath,amssymb,amsthm}
\usepackage[super]{nth}
\usepackage{algorithmic}
\usepackage[ruled,linesnumbered,noresetcount,vlined]{algorithm2e}
\usepackage{float}
\usepackage{bm}
\usepackage{framed}
\usepackage{arydshln} 

\ifCLASSINFOpdf
   \usepackage[pdftex]{graphicx}
\else

   \usepackage[dvips]{graphicx}
\fi

\ifCLASSOPTIONcompsoc
 \usepackage[caption=false,font=normalsize,labelfont=sf,textfont=sf]{subfig}
\else
 \usepackage[caption=false,font=footnotesize]{subfig}
\fi

\SetKwInput{KwInput}{Input}                
\SetKwInput{KwOutput}{Output}
\renewcommand{\baselinestretch}{1}

\floatstyle{ruled}
\newfloat{model}{thp}{lop}
\floatname{model}{Model}

\usepackage{enumitem}
\newlist{abbrv}{itemize}{1}
\setlist[abbrv,1]{label=,labelwidth=0.5in,align=parleft,itemsep=0.1\baselineskip,leftmargin=!}

\usepackage{csvsimple}
\usepackage{orcidlink}

\usepackage{tikz}
\usepackage{xcolor}
\usetikzlibrary{shapes.geometric, arrows.meta, positioning}

\tikzstyle{block} = [rectangle, rounded corners, minimum width=2.5cm, minimum height=1cm, text centered, draw=black, fill=gray!30]
\tikzstyle{arrow} = [thick,->,>=stealth]

\DeclareMathOperator*{\minimize}{minimize}

\begin{document}

\title{Improving the Accuracy of DC Optimal Power Flow Formulations via Parameter Optimization}

\author{\IEEEauthorblockN{Babak Taheri \orcidlink{0000-0002-7870-3465} and Daniel K. Molzahn \orcidlink{0000-0003-0583-5376}}
\thanks{\noindent Babak Taheri was with the School of Electrical and Computer Engineering, Georgia Institute of Technology, Atlanta, GA 30332 USA. He is now with Hitachi Energy Research, Raleigh, NC 27606 USA (email: babak.taheri@hitachienergy.com).

Daniel K. Molzahn is with the School of Electrical and Computer Engineering, Georgia Institute of Technology, Atlanta, GA 30332 USA (email: molzahn@gatech.edu). Support from NSF award \#2145564.}

}

\maketitle

\begin{abstract}

DC Optimal Power Flow (DC-OPF) problems optimize the generators' active power setpoints while satisfying constraints based on the DC power flow linearization. The computational tractability advantages of DC-OPF problems come at the expense of inaccuracies relative to AC Optimal Power Flow (AC-OPF) problems that accurately model the nonlinear steady-state behavior of power grids. This paper proposes an algorithm that significantly improves the accuracy of the generators' active power setpoints from DC-OPF problems with respect to the corresponding AC-OPF problems over a specified range of operating conditions. Using sensitivity information in a machine learning-inspired methodology, this algorithm tunes coefficient and bias parameters in the DC power flow approximation to improve the accuracy of the resulting DC-OPF solutions. Employing the Truncated Newton Conjugate-Gradient (TNC) method, a Quasi-Newton optimization technique, this parameter tuning occurs during an offline training phase, with the resulting parameters then used in online computations. Numerical results underscore the algorithm's efficacy with accuracy improvements in squared two-norm and $\infty$-norm losses of up to $90\%$ and $79\%$, respectively, relative to traditional DC-OPF formulations.

\end{abstract}


\section{Introduction}

Optimal Power Flow (OPF) is a fundamental tool for power system design and operation. OPF problems optimize system performance \textcolor{black}{(e.g., minimizing generation costs or reducing voltage deviations)}
while satisfying both equality constraints from a power flow model and inequality constraints from operational limits. Along with economic implications~\cite{cain2012history}, OPF problems are central to many algorithms for addressing uncertainty~\cite{roald2022review}, expansion planning~\cite{majidi2016stochastic}, improving stability margins~\cite{abhyankar2017tsopfIandII}, etc.

The AC power flow equations accurately model the steady-state behavior of power systems by relating the complex power injections with the voltage phasors. This motivates the solution of AC Optimal Power Flow (AC-OPF) problems that incorporate the AC power flow equations. However, AC-OPF problems are computationally challenging, with these problems being non-convex and NP-hard in the \textcolor{black}{worst case~\cite{bienstock2019strong,molzahn2017,bukhsh2013}}. 
These challenges are exacerbated for extensions of OPF problems that consider discrete characteristics, such as those related to transmission switching~\cite{barrows2014,vanAckooij2018} \textcolor{black}{ (where inaccuracies in DC model parameters can significantly impact the quality of switching decisions~\cite{taheri2024ac})} and unit commitment~\cite{padhy2004}, as well as uncertainties due to variable renewable energy generation and flexible loads~\cite{roald2022review}. 

OPF algorithms have been a subject of research for many decades~\cite{Momoh1999-1and2, aravena_molzahn_zhang_et_al-go_intro}.
Most existing algorithms seek a local optimum due to the non-convex nature of OPF problems.
Machine learning (ML) tools can address OPF challenges, offering computational speedups but often lacking interpretability\footnote{\textcolor{black}{In this context, interpretability refers to the ability to directly relate a model’s parameters and variables to the physical principles of the power system, a feature often lacking in more complex, black-box approaches.
}} and consistency with physical intuition~\cite{Duchesne2020}.

Due to these challenges, engineers often employ power flow approximations to obtain more tractable convex OPF formulations~\cite{molzahn2019}. Using a common linear approximation of the AC power flow equations known as the DC power flow~\cite{stott2009dc} yields the DC Optimal Power Flow (DC-OPF) problem, which is widely used for both short- and long-term planning purposes. While other advanced linearizations of the AC power flow equations exist~\cite{coffrin2016linear,molzahn2019}, our focus in this paper is on improving the accuracy of the \textit{DC-OPF formulation} due to its widespread usage.

\textcolor{black}{The DC power flow is derived under a set of standard assumptions: a flat voltage profile (magnitudes near 1.0 p.u.), small voltage angle differences, and negligible network losses and reactive power flows. These assumptions render the model computationally tractable and are often reasonably accurate for high-voltage transmission systems. However, approximation errors are introduced when these assumptions do not hold. While heuristics exist to mitigate some inaccuracies, such as embedding estimated power losses within the load data, such approaches may not be robust across diverse operating conditions~\cite{stott2009dc}.}
\textcolor{black}{Inaccuracies from the DC power flow assumptions are particularly pronounced when reactive power flows and voltage variations play a significant role, such as in heavily loaded systems~\cite{taylor1994power} or systems with high penetrations of renewable generators~\cite{milano2013overview}}. In such cases, the DC-OPF problem may not adequately capture the complex behavior of the power system, leading to suboptimal or infeasible decisions in operational and planning contexts~\cite{barrows2014,dvijotham_molzahn-cdc2016,baker2021solutions, Khodaei2010}.

\textcolor{black}{
While the standard DC-OPF formulation is widely used~\cite{stott2009dc}, its accuracy limitations have spurred the development of various alternatives~\cite{molzahn2019}. These can be broadly categorized into those that propose new model structures, those that tune parameters of existing models, and those that are purely based on machine learning techniques.
A significant body of work focuses on creating more accurate linear approximations of the AC power flow equations, e.g.,~\cite{Coffrin2014}. The survey in~\cite{molzahn2019} catalogs a wide array of these relaxations and approximations. While potentially more accurate, these models typically introduce new mathematical structures that pose significant practical implementation challenges, as the conventional DC-OPF structure is deeply embedded in numerous industry-grade and academic tools for complex tasks like unit commitment and transmission planning. Altering these problems' fundamental structure would require a significant overhaul of these established tools.
Another research direction, more aligned with our philosophy, focuses on enhancing existing models by tuning their coefficients and bias parameters. This includes optimizing the linearization point to account for supply and demand uncertainty~\cite{Hohmann2019}. This paradigm of enhancing classical, interpretable models has also proven effective for other approximations, such as the LinDistFlow model for distribution systems~\cite{Taheri2024LinDistFlow}. Our prior work used a single-level optimization to tune DC power flow parameters to better match AC power flow results~\cite{taheri2023optimizing}. While effective at improving the approximation's accuracy, it did not directly optimize for the end-goal of improving the OPF solution itself.
More recently, black-box machine learning models have been developed to predict AC-OPF solutions with high fidelity~\cite{Donti2021}, though often at the cost of interpretability.} 

\textcolor{black}{Our work, in contrast, occupies a middle ground. We leverage techniques from the machine learning domain, specifically differentiable optimization layers~\cite{Amos2017, Agrawal2019}, but our goal is to enhance a classical, interpretable model.
This paper distinguishes itself by intentionally preserving the standard DC-OPF formulation to ensure immediate compatibility with existing tools. Furthermore, we introduce a bilevel optimization framework that represents a significant conceptual advance. Rather than simply matching intermediate power flow values, our method tunes the DC-OPF parameters to directly minimize the error in the final optimal generator setpoints. This solution-focused approach carves a third path between pure machine learning and classical, fixed-parameter models, enhancing accuracy while preserving the interpretability and compatibility of the standard DC-OPF model.%
}

The DC power flow approximation is foundational to the DC-OPF problem. This approximation relates active power flow between buses $i$ and $j$, denoted $p_{ij}$, to the phase angle difference $\theta_i - \theta_j$ through a proportionality coefficient $b_{ij}$: $p_{ij} = b_{ij} (\theta_i - \theta_j)$. While typically derived from line parameters such as resistance $r_{ij}$ and reactance $x_{ij}$ via various heuristics, tailoring these parameter values to specific system and operating range characteristics can yield significant accuracy improvements. Common heuristics choose $b_{ij} = 1/x_{ij}$ or $b_{ij} = -\Im(1/(r_{ij} + \mathrm{j}x_{ij}))$, where $\Im(\cdot)$ extracts the imaginary part of a complex argument. When resistance is nonzero ($r_{ij} \neq 0$), these choices for $b_{ij}$ produce slightly different results, impacting the accuracy of the DC power flow approximation. Selecting the optimal $b_{ij}$ for a given application is not straightforward and depends on system characteristics and objectives~\cite{stott2009dc}. Additionally, the introduction of bias parameters allows further flexibility by adjusting power injections and flows to account for shunts, HVDC infeeds, phase shifts, and line losses~\cite{stott2009dc}. Two variants of the DC power flow—cold-start and hot-start—differ in their reliance on prior information when determining the $b_{ij}$ coefficients and bias parameters.

Our previous work in~\cite{taheri2023optimizing} introduced an algorithm that optimizes both the coefficient and bias parameters to improve DC power flow accuracy across a range of operating conditions. This optimization selects the DC power flow parameter values which minimize the mismatch between the line flows from the DC and AC power flow models over a training set of scenarios. We emphasize that the standard DC power flow parameters (e.g., $b_{ij} = 1/x_{ij}$) are themselves derived from simplifying assumptions that may not hold in practice. \textcolor{black}{
The goal of this parameter optimization is not to redefine the physical reactance or susceptance of a line. Rather, it is to find an effective set of parameter values that causes the DC power flow approximation to better match the results of the full AC power flow model. As discussed above, the standard DC power flow parameters are themselves derived from simplifying assumptions; we seek a different set of parameters that better serves the purpose of more accurately solving an optimization problem across a range of operating conditions.
}

While the previous paper~\cite{taheri2023optimizing} focused on optimizing the coefficients and biases within the DC power flow approximation itself, specialization of the parameters to specific applications can yield even larger accuracy advantages. 
Accordingly, this paper optimizes parameters for the DC power flow approximation in order to improve accuracy of DC-OPF solutions relative to AC-OPF solutions (as opposed to the focus on DC power flow versus AC power flow in our prior work). Specifically, we seek the DC power flow parameters that best reduce the discrepancies of the generator active power setpoints associated with the DC-OPF solutions relative to those from the AC-OPF solutions across a range of load demands.  
\textcolor{black}{This represents a significant conceptual and technical advancement that necessitates reformulating as a bilevel optimization problem that differs substantially from the single-level problem in~\cite{taheri2023optimizing}. The upper level seeks to tune the DC-OPF parameters, while the lower level solves the convex DC-OPF problem itself. This enhanced structure introduces new challenges, as the loss function must consider the optimal operational decisions from the lower-level optimization.}

Leveraging differentiable optimization techniques developed for deep learning architectures~\cite{Amos2017,Agrawal2019},  this paper proposes a new algorithm for adaptively selecting coefficients $b_{ij}$ and bias parameters within the specific context of DC-OPF problems. The algorithm performs offline fine-tuning of $b_{ij}$ and bias parameters using \textcolor{black}{a training dataset of AC-OPF solutions under various representative operating conditions}. We differentiate the lower-level problem to compute the sensitivities of the DC-OPF solutions with respect to these parameters. This enables us to iteratively adjust the parameter values during an offline phase by using a gradient-based optimization method, specifically the Truncated Newton Conjugate-Gradient (TNC) algorithm, to minimize the difference between DC-OPF's generator active power setpoints and their AC-OPF counterparts across the training dataset. The resulting optimized parameter values are then applied to improve the accuracy of the DC-OPF problems during online computations.

To summarize, the main contributions of this paper are:
\begin{enumerate}
    \item \textcolor{black}{Formulating the task of tuning DC-OPF parameters as a bilevel optimization problem, where the goal is to minimize the error in the optimal generator setpoints relative to full AC-OPF solutions.}

    \item Introducing an optimization algorithm designed to adaptively choose the coefficients $b_{ij}$ and bias parameters of the DC power flow inside the DC-OPF problem. These parameters are instrumental in modeling line losses and incorporating factors such as shunts, HVDC infeeds, phase shifts, and line losses.

     \item Computing gradients for the DC-OPF solution with respect to both coefficient and bias parameters to obtain sensitivities of the generator setpoints as well as our defined loss function, \textcolor{black}{enabled by modern differentiable programming techniques}.

    \item Utilizing a quasi-Newton method (TNC) to scale our proposed algorithm to large power systems.
    
    \item Providing numerical results that demonstrate the superior accuracy of our proposed algorithm over a range of operating conditions.
\end{enumerate}

The rest of the paper is structured as follows. Section~\ref{sec:OPF} reviews the AC-OPF and DC-OPF problems. Section~\ref{sec:Proposed_Algorithm} presents our algorithm for optimizing parameters of DC-OPF problems. Section~\ref{sec:Numerical experiments} provides numerical results to demonstrate the algorithm's performance. Section~\ref{sec:conclusion} concludes the paper.

\section{Optimal Power Flow Problem}
\label{sec:OPF}

To introduce notation and necessary background material, this section describes the AC-OPF formulation and the DC-OPF linear approximation that simplifies the AC-OPF problem to improve tractability at the cost of accuracy. Buses, lines and generators in the network are denoted by the sets $\mathcal{N}$, $\mathcal{E}$ and $\mathcal{G}$, respectively. \textcolor{black}{Here, $\mathcal{E}$ represents the set of individual branches; we note that multiple parallel branches may exist between two buses, and our formulation treats each one distinctly with its own parameters and flow limits.} Each bus~$i\in\mathcal{N}$ has a voltage phasor $V_i$ with phase angle $\theta_{i}$, along with a shunt admittance $Y_i^S$, complex power demand $s^{\text{d}}_i =\text{p}^{\text{d}}_i +\mathrm{j}\text{q}^{\text{d}}_i $, and generated complex power $s^{\text{g}}_i = p^{\text{g}}_i + \mathrm{j}q^{\text{g}}_i$, where $\mathrm{j}$ is the imaginary unit $\mathrm{j} =\sqrt{-1}$. Buses without generators are modeled as having upper and lower generation limits of zero. The complex power flows entering each end of line $(j,k)\in\mathcal{E}$ are represented by $S_{jk}$ and $S_{kj}$. Series admittance parameters for each line $(j,k)$ are $Y_{jk}$ and $Y_{kj}$ and the line's shunt admittance is $Y^{c}_{jk}$. The real and imaginary components of a complex number are indicated by $\Re(\,\cdot\,)$ and $\Im(\,\cdot\,)$, respectively, $\|\cdot\|_{\infty}$ denotes the $L_{\infty}$-norm, and $\|\cdot\|_{2}$ denotes the $L_{2}$-norm. The complex conjugate and the transpose of a matrix are denoted by $(\,\cdot\,)^{\star}$ and $(\,\cdot\,)^{\top}$, respectively, and $\angle (\,\cdot\,)$ is the angle of a complex argument. Upper and lower limits are represented by $(\overline{\,\cdot\,})$ and $(\underline{\,\cdot\,})$. 

\subsection{The AC-OPF Problem}
Model~\ref{model:ACOPF} presents the AC-OPF formulation. The objective~\eqref{eq:ACOPF:objective} minimizes the total generation cost by summing the costs for each generator $i\in\mathcal{N}$, modeled as quadratic functions (with coefficients $c_{2i}$, $c_{1i}$, and $c_{0i}$) of the real part of the complex power generation \( s_i^{\text{g}} \). \textcolor{black}{Following typical modeling practices, we assume that $c_{2i} \geq 0$ such that the generator cost function is convex.}
Constraint~\eqref{eq:ACOPF:voltage_bounds} ensures that all voltage magnitudes are within specified limits. Constraint~\eqref{eq:ACOPF:dispatch_bounds} bounds the complex power outputs of each generator within its feasible operating range, where inequalities on complex quantities are interpreted as bounds on the real and imaginary parts. Constraint~\eqref{eq:ACOPF:thermal_limits} enforces the line flow limits.

\label{sec:OPF:AC}

    \begin{model}[!t]
        \caption{The AC-OPF Problem}
        \label{model:ACOPF}
        \begin{subequations}
        \label{eq:ACOPF}
        
        \begin{align}
            \minimize_{s^{\text{g}}_{i}, V_{i}} \quad 
            & \sum_{i \in \mathcal{N}} c_{2i} (\Re(s^{\text{g}}_{i}))^{2 }+ c_{1i}\Re(s^{\text{g}}_{i}) + c_{0i} \label{eq:ACOPF:objective}\\
            \textrm{s.t.} \quad  &(\forall i\in\mathcal{N}, ~\forall (j,k)\in\mathcal{E}) \nonumber\\
                & \underline{\text{V}}_{i} \leq |V_{i}| \leq \overline{\text{V}}_{i} 
                &&  \label{eq:ACOPF:voltage_bounds}\\
                & \underline{s}^{\text{g}}_{i} \leq s^{\text{g}}_{i} \leq \overline{s}^{\text{g}}_{i}
                &&  \label{eq:ACOPF:dispatch_bounds}\\
                &  \left|S_{jk}\right| \leq \overline{S}_{jk}, ~\left|S_{kj}\right| \leq \overline{S}_{jk}
                &&  \label{eq:ACOPF:thermal_limits} \\
                & s^{\text{g}}_{i} - s^{\text{d}}_{i} - (Y^{s}_{i})^{\star} |V_{i}|^{2} =\!\!\sum_{(i, j) \in \mathcal{E}}\!\! S_{i j} + \!\!\sum_{(k, i) \in \mathcal{E}}\!\! S_{ik}
                &&  \label{eq:ACOPF:kirchhoff} \\
            & S_{jk} = (Y_{jk} + Y_{jk}^{c})^{\star} V_{j} V_{j}^{\star} - Y_{jk}^{\star} V_{j} V_{k}^{\star}
                &&  \label{eq:ACOPF:ohm_fr}\\
            & S_{kj} = (Y_{kj} + Y_{kj}^{c})^{\star} V_{k} V_{k}^{\star} - Y_{kj}^{\star} V_{k} V_{j}^{\star}
                &&  \label{eq:ACOPF:ohm_to}\\    
            &\theta_{ref}=0 &&  \label{eq:ACOPF:ref}\\  
            & - \overline{\theta}_{jk} \leq \angle (V_{j} V_{k}^{\star}) \leq\overline{\theta}_{jk}  &&  \label{eq:ACOPF:phase_angle_limits}\\
            \textbf{variables}: & ~s^{\text{g}}_{i}~(\forall i \in \mathcal{N}), V_{i}~(\forall i \in \mathcal{N})  \nonumber
        \end{align}
        \end{subequations}
    \end{model}

Constraint~\eqref{eq:ACOPF:kirchhoff} ensures power balance at each bus and \eqref{eq:ACOPF:ohm_fr} and~\eqref{eq:ACOPF:ohm_to} represent the power flows on each line.
The reference angle is set to zero in \eqref{eq:ACOPF:ref}. Finally,~\eqref{eq:ACOPF:phase_angle_limits} limits the phase angle difference across each line.
The AC-OPF problem is inherently nonlinear and non-convex due to the quadratic nature of the power flow equations. Using nonlinear optimization techniques such as interior-point methods, AC-OPF algorithms frequently find high-quality solutions~\cite{aravena_molzahn_zhang_et_al-go_intro, Momoh1999-1and2, gopinath2022}.

\subsection{The DC-OPF Problem}
\label{sec:OPF:DC}
The DC-OPF problem is a linear approximation of the AC-OPF problem that assumes uniform voltage magnitudes (typically normalized to one per-unit), small voltage angle differences, negligible losses, and ignores reactive power~\cite{stott2009dc}. These assumptions are most relevant to transmission systems. The DC-OPF model underlies various applications within electricity markets and is widely used in unit commitment, transmission network expansion planning, and optimal transmission switching problems~\cite{vanAckooij2018,knueven2020mixed, fisher2008optimal, hedman2008optimal, guo2022tightening, Khodaei2010}.
    \begin{model}[!th]
        \caption{The DC-OPF Problem}
        \label{model:DCOPF}
        \begin{subequations}
        \label{eq:DCOPF}
        \normalsize
        \begin{align}
           \minimize_{p^{\text{g}}_{i}, \theta_{i}} \quad 
            & \sum_{i \in \mathcal{N}} c_{2i} (p^{\text{g}}_{i})^{2 }+ c_{1i}p^{\text{g}}_{i} + c_{0i}  \label{eq:DCOPF:obj}\\
            \text{s.t.} \quad &(\forall i\in\mathcal{N}, ~\forall (j,k)\in\mathcal{E}) \nonumber\\
            & \underline{\text{p}}^{\text{g}}_{i} \leq p^{\text{g}}_{i} \leq \overline{\text{p}}^{\text{g}}_{i} 
                   \label{eq:DCOPF:bounds:pg}\\
            & |P_{jk}| \leq \bar{S}_{jk}
                    \label{eq:DCOPF:bounds:pf}\\       
                & p^{\text{g}}_{i} - \text{p}^{\text{d}}_{i} - \gamma_{i} = \sum_{(i,j) \in \mathcal{E}} P_{ij} - \sum_{(k,i) \in \mathcal{E}} P_{ki}   
                   \label{eq:DCOPF:power_balance}\\
                & P_{jk} = b_{jk} (\theta_{j} - \theta_{k}) + \rho_{jk}
                   \label{eq:DCOPF:ohm}\\
                   & \theta_{ref} = 0 \label{eq:DCOPF:ref}\\
                \textbf{variables}: & ~p^{\text{g}}_{i}~(\forall i \in \mathcal{N}), \theta_{i}~(\forall i \in \mathcal{N})  \nonumber
        \end{align}
        \end{subequations}
    \end{model}

Model~\ref{model:DCOPF} presents the DC-OPF problem. Analogous to~\eqref{eq:ACOPF}, \eqref{eq:DCOPF:obj} minimizes the total generation cost,~\eqref{eq:DCOPF:bounds:pg} bounds the generators' active power outputs,~\eqref{eq:DCOPF:bounds:pf} enforces line flow limits, and~\eqref{eq:DCOPF:power_balance} balances active power at each bus.
Constraint~\eqref{eq:DCOPF:ohm} linearizes the line flow expressions by relating branch power flows to voltage angle differences according to the DC power flow approximation, thus offering a simplified counterpart to~\eqref{eq:ACOPF:ohm_fr} and~\eqref{eq:ACOPF:ohm_to}. Finally,~\eqref{eq:DCOPF:ref} sets the reference angle.

\textcolor{black}{Note that \eqref{eq:DCOPF:power_balance} and \eqref{eq:DCOPF:ohm} include parameters $\boldsymbol{\gamma}$ and $\boldsymbol{\rho}$ to improve accuracy. Following extended DC power flow models~\cite{stott2009dc}, these parameters have distinct roles. The flow bias $\rho_{ij}$ in~\eqref{eq:DCOPF:ohm} can account for factors like phase-shifter effects and active power losses on the line $(i,j)$. The injection bias $\gamma_i$ in~\eqref{eq:DCOPF:power_balance} can account for losses from shunts, HVDC infeeds, or aggregate the effects of losses from all connected lines.}

Model~\ref{model:DCOPF_matrixformat} equivalently rewrites Model~\ref{model:DCOPF} in a matrix form, where $\mathbf{A}$ is the $|\mathcal{E}| \times |\mathcal{N}|$ branch-bus incidence matrix and $\text{diag}(\mathbf{b})$ is the diagonal matrix with $\mathbf{b}$ on the diagonal. 

    \begin{model}[!t]
        \caption{The DC-OPF-M Problem}
        \label{model:DCOPF_matrixformat}
        \begin{subequations}
        \label{eq:DCOPF_matrixformat}
        \normalsize
        \begin{align}
           \minimize_{\mathbf{p}^{\text{g}} , \boldsymbol{\theta}} \quad 
            & {\mathbf{p}^{\text{g}}}^{\top}\text{diag}(\mathbf{c}_{2}) \mathbf{p}^{\text{g}}+ \mathbf{c}^{\top}_{1}\mathbf{p}^{\text{g}} + \mathbf{1}^{\top}\mathbf{c}_{0} \label{eq:DCOPF:obj_matrixformat}\\
            \text{s.t.} \quad \nonumber\\
             & \underline{\textbf{p}}^{\text{g}} \leq \mathbf{p}^{\text{g}} \leq \overline{\textbf{p}}^{\text{g}}
                   \label{eq:DCOPF:bounds:pg_matrixformat}\\
                &  \mathbf{p}^{\text{g}} - \textbf{p}^{\text{d}} - \boldsymbol{\gamma} = \mathbf{A}^{\top} \Big(\text{diag}(\mathbf{b}) \mathbf{A} \boldsymbol{\theta}  +\boldsymbol{\rho}\Big) 
                   \label{eq:DCOPF:power_balance_matrixformat}\\
                & |\text{diag}(\mathbf{b})\mathbf{A} \boldsymbol{\theta} + \boldsymbol{\rho}| \leq \bar{\mathbf{S}}
                    \label{eq:DCOPF:bounds:pf_matrixformat}\\ 
                    & \theta_{ref} = 0\\
                \textbf{variables}: & ~\mathbf{p}^{\text{g}} , \boldsymbol{\theta}  \nonumber
        \end{align}
        \end{subequations}
    \end{model}

The DC-OPF problem's accuracy is dictated by the parameters $\mathbf{b}$, $\boldsymbol{\gamma}$, and $\boldsymbol{\rho}$. Traditionally, there are two DC-OPF versions, cold-start and hot-start, that rely on different amounts of prior information to select these parameters.

\subsubsection{Cold-start DC-OPF}
In this version, the coefficient and bias parameters are selected without relying on a nominal operating point. For instance, the coefficient values $b_{ij}$ can be selected as the imaginary part of the line susceptance:
\begin{equation}
\label{eq:cold-start}
b^{cold}_{ij} = \Im\left(\frac{-1}{r_{ij} + \mathrm{j} x_{ij}}\right).
\end{equation}
We note that the alternative selection of $b_{ij}^{cold} = 1/x_{ij}$ yields very similar numerical results to the choice of $b_{ij}^{cold}$ in~\eqref{eq:cold-start}, so we only present results for the latter for the sake of brevity.

The cold-start version assigns values of zero to the bias parameters $\boldsymbol{\gamma}$ and $\boldsymbol{\rho}$. The simplicity provided by the cold-start DC-OPF problem comes at the cost of accuracy in settings where information on a nominal operating point is available.

\subsubsection{Hot-start DC-OPF}
Hot-start DC-OPF problems leverage information from a nominal AC power flow solution to improve accuracy. There are many hot-start DC-OPF variants~\cite{stott2009dc}; e.g., the ``localized loss modeling'' variant in~\cite{stott2009dc} selects the coefficient and bias parameters as:
\begin{subequations}
    \begin{align}
        b^{\text{hot}}_{ij} &= b_{ij}\, |V^{\bullet}_{i}|\, |V^{\bullet}_{j}| \frac{\sin(\theta^{\bullet}_{i} - \theta^{\bullet}_{j})}{\theta^{\bullet}_{i} - \theta^{\bullet}_{j}}, \label{eq:sub1-hot-start} \\
        \gamma^{hot}_{i} &= \sum_{(i,j) \in \mathcal{E}} \Re(Y_{ij}) |V^{\bullet}_{i}|(|V^{\bullet}_{i}|-|V^{\bullet}_{j}| \cos(\theta^{\bullet}_{i}-\theta^{\bullet}_{j})), \label{eq:sub2-hot-start}\\
        \rho^{hot}_{ij}&= \Re(Y_{ij}) |V^{\bullet}_{i}|(|V^{\bullet}_{i}|-|V^{\bullet}_{j}| \cos(\theta^{\bullet}_{i}-\theta^{\bullet}_{j})), \label{eq:sub3-hot-start}
    \end{align}%
    \label{eq:hot-start}%
\end{subequations}%
where the notation $(\,\cdot\,)^{\bullet}$ signifies quantities associated with a nominal AC power flow solution, \textcolor{black}{and \( b_{ij} \) on the right-hand side of~\eqref{eq:sub1-hot-start} is typically the cold-start version defined in~\eqref{eq:cold-start}}. The bias $\boldsymbol{\gamma}^{hot}$ aims to capture the effects of branch losses on phase angles and $\boldsymbol{\rho}^{hot}$ incorporates the branch losses in the line flow expressions.
\textcolor{black}{
A notable limitation of the hot-start approach is its reliance on a single, pre-determined nominal operating point, denoted by the \((\bullet)\) superscript in its defining equations. The accuracy of the hot-start parameters is therefore highest near that specific point and can degrade significantly as conditions deviate. This presents a challenge when a DC power flow approximation is needed for a wide range of operating conditions, as selecting hot-start parameters based on a single nominal operating point may yield insufficient accuracy.
}

\subsubsection{Optimized Parameters for DC Power Flow (DCPF)~\cite{taheri2023optimizing}}

Our recent work in~\cite{taheri2023optimizing} proposes an alternative for selecting the $\mathbf{b}$, $\boldsymbol{\gamma}$, and $\boldsymbol{\rho}$ parameters by solving a single-level optimization problem. Providing a conceptual foundation for the approach proposed in this paper, the algorithm in~\cite{taheri2023optimizing} also optimizes these parameters over a range of operation. However, there is a key distinction between the algorithm in~\cite{taheri2023optimizing} and the algorithm in this paper: the approach in~\cite{taheri2023optimizing} does not explicitly consider the impacts of the optimized parameters on the solution to a DC-OPF problem, but rather seeks to balance the DC power flow's accuracy with respect to all operating points within a specified range. Conversely, our proposed algorithm in this paper leverages machine learning techniques for training models with differentiable optimization layers to tailor the DC power flow parameter values so that the DC \emph{optimal} power flow solutions align with the AC-OPF solutions. As shown numerically in Section~\ref{sec:Numerical experiments}, this enables substantial accuracy improvements relative to~\cite{taheri2023optimizing} in settings where the optimized DC power flow parameters are used in optimization problems like DC-OPF.
\textcolor{black}{This represents a significant conceptual and technical advancement over the single-level problem in~\cite{taheri2023optimizing}. While the prior work focused on tuning parameters to match the line flow behaviors of the DC and AC power flow \textit{equations}, the present work introduces a bilevel optimization framework. The objective is no longer to fit the intermediate equations, but rather to tune the parameters such that the final \textit{optimal solution} of the DC-OPF problem aligns with that of the AC-OPF problem. This solution-focused methodology is a key contribution of this paper.%
}

We next present our algorithm for optimizing the coefficient ($\mathbf{b}$) and bias ($\boldsymbol{\gamma}$ and $\boldsymbol{\rho}$) parameters. Our algorithm seeks to match the generators' active power outputs from the DC-OPF solution to corresponding values from the AC-OPF solution.

\section{Parameter Optimization Algorithm}
\label{sec:Proposed_Algorithm}

\begin{figure}[!t]
\centerline{\includegraphics[width=0.48\textwidth]{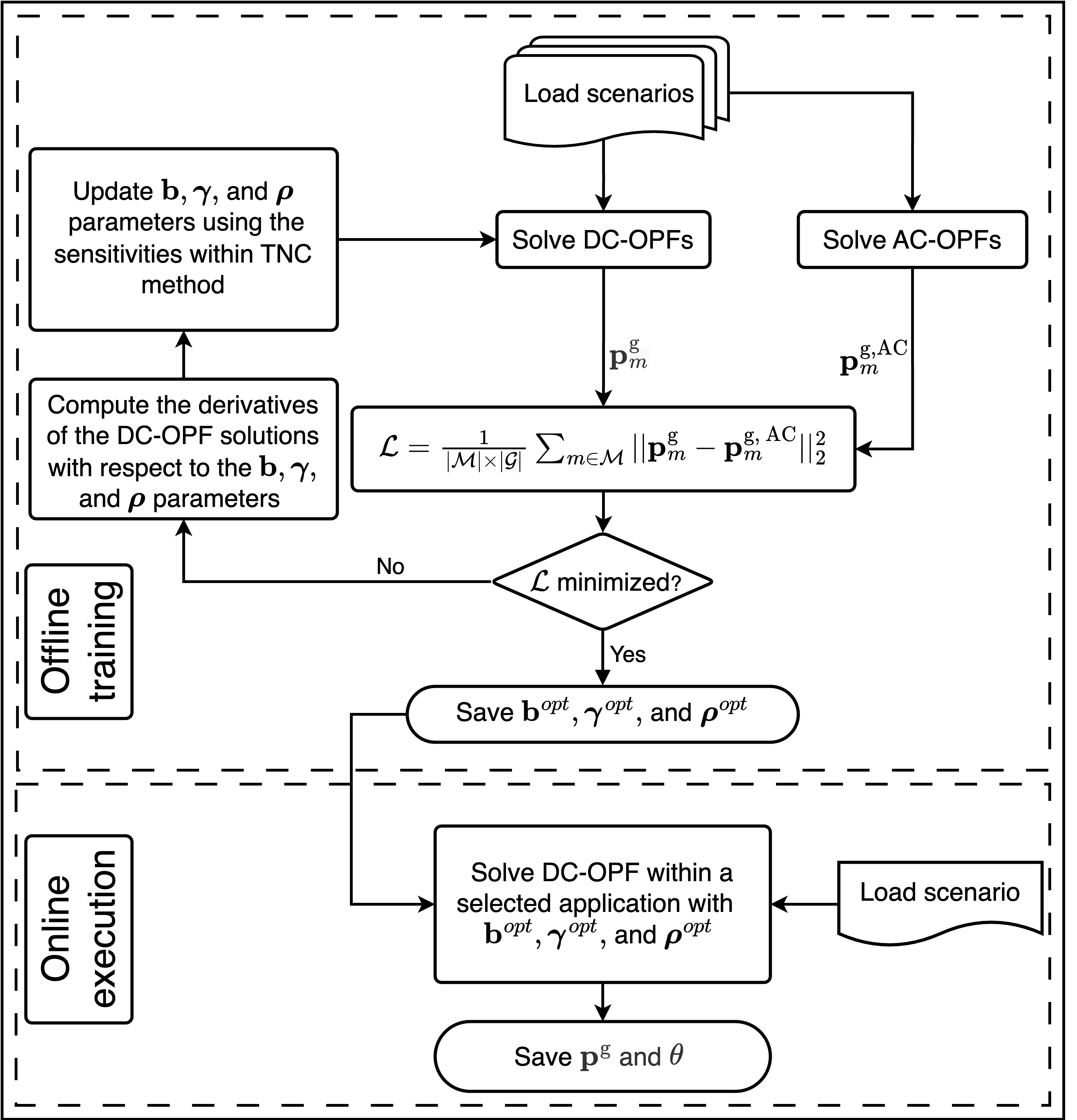}}
\caption{Flowchart of the proposed algorithm.}
\label{fig:flowchart}
\vspace{-1em}
\end{figure}

Our parameter optimization algorithm, depicted in Fig.~\ref{fig:flowchart}, has offline and online phases. The offline phase computes DC-OPF parameters (\(\mathbf{b}\), \(\boldsymbol{\gamma}\), and \(\boldsymbol{\rho}\)) to achieve a high-fidelity approximation of the AC-OPF solution over a range of operational conditions. The online phase employs these parameters to efficiently solve DC-OPF problems in various applications. While our numerical results focus on DC-OPF problems, the optimized DC-OPF parameters are applicable to many computationally intensive optimization problems that incorporate DC power flow models, such as unit commitment~\cite{vanAckooij2018,knueven2020mixed}, expansion planning~\cite{majidi2016stochastic}, and grid reliability enhancement~\cite{Liu2010}.

Our algorithm minimizes a loss function based on the discrepancy between DC-OPF and AC-OPF generator active power setpoints over various load scenarios.
To perform this minimization, we compute the loss function's gradients with respect to the parameters $\mathbf{b}$, $\boldsymbol{\gamma}$, and $\boldsymbol{\rho}$ for use in a nonlinear optimization method.
Based on our numerical tests comparing various optimization methods, TNC's approach to approximating the Hessian matrix from gradient information best achieved scalable and effective parameter optimization.

\subsection{Loss Function}
\label{subsec:Loss Function}
We first introduce a loss function $\mathcal{L}$ that sums the squared two-norm discrepancies between the generators' active power setpoints from the AC-OPF problems ($\mathbf{p}_{m}^{\text{g,AC}}$) and the DC-OPF problems ($\mathbf{p}_{m}^{\text{g}}$) across a specified set of load scenarios $\mathcal{M}$:
\begin{align}
\mathcal{L}(\mathbf{b}, \boldsymbol{\gamma}, \boldsymbol{\rho}) &= \frac{1}{|\mathcal{G}|\times|\mathcal{M}|} \sum_{m \in \mathcal{M}} \|\mathbf{p}_{m}^{\text{g}}(\mathbf{b},\boldsymbol{\gamma},\boldsymbol{\rho}) - \mathbf{p}_{m}^{\text{g,AC}}\|^2_2,\label{eq:objective_function}%
\end{align}%
where the constant~$\frac{1}{|\mathcal{G}|\times|\mathcal{M}|}$ normalizes this function based on the number of generators and load scenarios. As shown in~\eqref{eq:objective_function}, $\mathbf{p}_m^{\text{g}}$ (and thus $\mathcal{L}(\mathbf{b}, \boldsymbol{\gamma},\boldsymbol{\rho})$) is a function of the coefficient parameters $\mathbf{b}$ and the bias parameters $\boldsymbol{\gamma}$ and $\boldsymbol{\rho}$.
\textcolor{black}{
The selection of the training scenarios in the set $\mathcal{M}$ is an important application-dependent step. Scenarios should be chosen to represent the operating conditions under which the optimized parameters will be used. In our experiments, we generate scenarios by scaling nominal loads based on a normal distribution; however, in practice, utilities could use historical data or forecast scenarios to create a training set tailored to their system.
}

We note that this two-norm loss formulation is typical in machine learning for its robustness and differentiability. Moreover, this loss function penalizes larger deviations more heavily, which is well aligned with typical applications where a small number of severe approximation errors would be more problematic than a large number of minor errors. One could instead use other norms without major conceptual changes.

With this loss function, the optimization problem to find the best coefficient and bias parameters is formulated as:
\begin{equation}
\minimize_{\mathbf{b},\boldsymbol{\gamma},\boldsymbol{\rho}} \quad \mathcal{L}(\mathbf{b},\boldsymbol{\gamma},\boldsymbol{\rho}).
    \label{eq:optimization problem}
\end{equation}
Using the sensitivity analysis that we describe in the next subsection, we adjust the coefficient and bias parameters ($\mathbf{b}$, $\boldsymbol{\gamma}$, and $\boldsymbol{\rho}$) to optimize the DC-OPF problem's parameters. As detailed in Model~\ref{model:O-DCOPF_matrixformat}, this task takes the form of a bilevel optimization problem, with the upper-level fine-tuning parameters based on outcomes from lower-level DC-OPF problems under a variety of operational scenarios.

    \begin{model}[!t]
        \caption{Optimizing DC-OPF Problem}
        \label{model:O-DCOPF_matrixformat}
        \begin{subequations}
        \label{eq:O-DCOPF_matrixformat}
        \normalsize
        \begin{align}
         \minimize_{\mathbf{b}, \boldsymbol{\gamma}, \boldsymbol{\rho}} \quad & \frac{1}{|\mathcal{G}|\times |\mathcal{M}|}\sum_{m \in \mathcal{M}} \|\mathbf{p}_{m}^{\text{g}}(\mathbf{b},\boldsymbol{\gamma},\boldsymbol{\rho}) - \mathbf{p}^{\text{g},\text{AC}}_{m}\|^2_2 \label{eq:O-DCOPF:main_objective}\\ %
     \text{s.t.}\quad \quad \quad &  \nonumber \\
             \minimize_{\mathbf{p}^{\text{g}}_{m} , \boldsymbol{\theta}_{m}} &\quad  {\mathbf{p}^{\text{g}}_{m}}^{\top}\text{diag}(\mathbf{c}_{2}) \mathbf{p}^{\text{g}}_{m}+ \mathbf{c}^{\top}_{1}\mathbf{p}^{\text{g}}_{m} + \mathbf{1}^{\top}\mathbf{c}_{0}  \label{eq:O-DCOPF:obj_opt}\\
            \text{s.t.}  & \nonumber\\
                & \mathbf{p}^{\text{g}}_{m} - \textbf{p}^{\text{d}}_{m} - \boldsymbol{\gamma} = \mathbf{A}^{\top} \Big(\text{diag}(\mathbf{b}) \mathbf{A} \boldsymbol{\theta}_{m} +\boldsymbol{\rho}\Big)   
                   \label{eq:O-DCOPF:power_balance_matrixformat}\\
                & \quad \quad |\text{diag}(\mathbf{b})\mathbf{A} \boldsymbol{\theta}_{m} + \boldsymbol{\rho}| \leq \bar{\mathbf{S}}
                    \label{eq:O-DCOPF:bounds:pf_matrixformat}\\
                & \quad \quad \underline{\textbf{p}}^{\text{g}} \leq \mathbf{p}^{\text{g}}_{m} \leq \overline{\textbf{p}}^{\text{g}}
                   \label{eq:O-DCOPF:bounds:pg_matrixformat}\\
                   &\quad \quad \theta_{ref} = 0 \label{eq:O-DCOPF:angle_ref}\\
                   &\quad  \forall m \in \mathcal{M} \nonumber
        \end{align}
        \end{subequations}
        
    \end{model}

\begin{algorithm}[!t]
\caption{Truncated Newton (TNC) Method}
\label{alg:truncated_newton_with_preconditioning}
\small
\DontPrintSemicolon

\KwInput{\hspace{0cm}$\mathbf{x}_0 = [\mathbf{b}^{\top}_{0}, \boldsymbol{\gamma}^{\top}_{0}, \boldsymbol{\rho}^{\top}_{0}]^{\top}$: Initial guess\\
\hspace{1cm}$\epsilon$: Tolerance for convergence \\
\hspace{1cm}$\textit{max\_iter}$: Maximum iterations\\
\hspace{1cm}$\mathcal{L}(\mathbf{x}_{k})$: Loss function\\
\hspace{1cm}$\nabla \mathcal{L}(\mathbf{x}_{k}) = \mathbf{g}$\\
\hspace{1cm}$\mathbf{M}$: Preconditioning matrix (often diagonal) \\
\hspace{1cm}$\mathbf{H}$: Hessian or its approximation function \\
\hspace{1cm}$\alpha_1$: Armijo condition constant (e.g., $10^{-4}$) \\
\hspace{1cm}$\alpha_2$: Curvature condition constant, between $\alpha_1$ and $1$ 
}

\KwOutput{
    Optimized parameters $\mathbf{x}^*$
}

Initialize $\mathbf{x}_k \leftarrow \mathbf{x}_0$, $k \leftarrow 0$\;

\While{$k \leq \textit{max\_iter}$ and $\|\nabla \mathcal{L}(\mathbf{x}_{k})\| > \epsilon$}
{
    $\mathbf{g} \gets \nabla \mathcal{L}(\mathbf{x}_{k})$\;
    $\mathbf{r} \gets -\mathbf{g}$\;
    $\mathbf{z} \gets \mathbf{M}^{-1} \mathbf{r}$\;
    $\mathbf{p} \gets \mathbf{z}$\;
    $\mathbf{s} \gets \mathbf{0}$\;
    $\rho_{old} \gets \mathbf{r}^{\top} \mathbf{z}$\;
    
    \While{$\|\mathbf{r}\| > \epsilon$}{
        $\mathbf{q} \gets \mathbf{H} \mathbf{p}$\;
        $\alpha \gets \frac{\rho_{old}}{\mathbf{p}^{\top} \mathbf{q}}$\;
        $\mathbf{s} \gets \mathbf{s} + \alpha \mathbf{p}$\;
        $\mathbf{r} \gets \mathbf{r} - \alpha \mathbf{q}$\;
        $\mathbf{z} \gets \mathbf{M}^{-1} \mathbf{r}$\;
        $\rho_{new} \gets \mathbf{r}^{\top} \mathbf{z}$\;
        $\eta \gets \frac{\rho_{new}}{\rho_{old}}$\;
        $\mathbf{p} \gets \mathbf{z} + \eta \mathbf{p}$\;
        $\rho_{old} \gets \rho_{new}$\;
    }
    
    \tcp*{Wolfe Line Search to determine $\beta$}
    $\beta \gets 1$\;
    \While{True}{
        \If{$\mathcal{L}(\mathbf{x}_{k} + \beta \mathbf{s}) \leq \mathcal{L}(\mathbf{x}_{k}) + \alpha_1 \beta \mathbf{g}^{\top} \mathbf{s}$ \textbf{and} $\|\nabla \mathcal{L}(\mathbf{x}_{k} + \beta \mathbf{s})^{\top} \mathbf{s}\| \leq \alpha_2 \|\mathbf{g}^{\top} \mathbf{s}\|$}{
            \textbf{break}\;
        }
        $\beta \gets \beta / 2$\; 
    }
    $\mathbf{x}_{k+1} \gets \mathbf{x}_{k} + \beta \mathbf{s}$\;
    $k \gets k + 1$\;
}

$\mathbf{x}^{*} \leftarrow \mathbf{x}_{k}$\;

\end{algorithm}

\subsection{Sensitivity Analysis of Coefficient and Bias Parameters}
\label{subsec:Sensitivity Analysis}

This subsection describes the sensitivity analyses crucial for our parameter optimization algorithm. We utilize a gradient-based method, TNC, which requires the gradient of the loss function with respect to $\mathbf{b}$, $\boldsymbol{\gamma}$, and $\boldsymbol{\rho}$. These gradients, denoted as $\mathbf{g}^{b}$, $\mathbf{g}^{\gamma}$, and $\mathbf{g}^{\rho}$, are computed as follows: 
\begin{subequations}
\label{eq:sensitivity}
\begin{align}
\label{eq:sensitivity_b}
 \mathbf{g}^{b}=\frac{\partial \mathcal{L}}{\partial \mathbf{b}}=  \frac{2}{|\mathcal{G}| \times |\mathcal{M}|}\sum_{m \in \mathcal{M}} \left. \Big(\frac{\partial \mathbf{p}^{\text{g}}}{\partial \mathbf{b}} \right|_{\mathbf{p}^{\text{g}}_{m}}\Big)^{\top} \Big(\mathbf{p}_{m}^{\text{g}} - \mathbf{p}_{m}^{\text{g},AC} \Big),
\end{align}
\begin{align}
\label{eq:sensitivity_gamma}
 \mathbf{g}^{\gamma}=\frac{\partial \mathcal{L}}{\partial \boldsymbol{\gamma}}=  \frac{2}{|\mathcal{G}|\times |\mathcal{M}|}\sum_{m \in \mathcal{M}} \left. \Big(\frac{\partial \mathbf{p}^{\text{g}}}{\partial \boldsymbol{\gamma}} \right|_ {\mathbf{p}^{\text{g}}_{m}}\Big)^{\top}\Big(\mathbf{p}_{m}^{\text{g}} - \mathbf{p}_{m}^{\text{g},AC} \Big),
\end{align}
\begin{align}
\label{eq:sensitivity_rho}
 \mathbf{g}^{\rho}=\frac{\partial \mathcal{L}}{\partial \boldsymbol{\rho}}=  \frac{2}{|\mathcal{G}|\times |\mathcal{M}|}\sum_{m \in \mathcal{M}} \left. \Big(\frac{\partial \mathbf{p}^{\text{g}}}{\partial \boldsymbol{\rho}} \right|_ {\mathbf{p}^{\text{g}}_{m}}\Big)^{\top}\Big(\mathbf{p}_{m}^{\text{g}} - \mathbf{p}_{m}^{\text{g},AC} \Big).
\end{align}
\end{subequations}

\textcolor{black}{The Jacobian matrices $\frac{\partial \mathbf{p}^{\text{g}}}{\partial \mathbf{b}}$, $\frac{\partial \mathbf{p}^{\text{g}}}{\partial \boldsymbol{\gamma}}$, and $\frac{\partial \mathbf{p}^{\text{g}}}{\partial \boldsymbol{\rho}}$ are central to the gradients in~\eqref{eq:sensitivity}. Since $\mathbf{p}^{\text{g}}$ is the solution to the convex DC-OPF optimization problem, these derivatives are found by implicitly differentiating the Karush-Kuhn-Tucker (KKT) optimality conditions of the DC-OPF problem~\cite{fiacco1990nonlinear, amos2017optnet}. This is a key step that makes our bilevel approach tractable.}

\textcolor{black}{To perform these implicit derivative computations, we leverage techniques developed for embedding convex optimization layers within deep learning architectures; specifically, we utilize \texttt{cvxpylayers}\cite{cvxpylayers2019}, which allows for automatic differentiation through convex optimization problems specified in \texttt{CVXPY}. This integration facilitates efficient and scalable gradient computations necessary for our optimization framework. Our implementation defines the DC-OPF problem using CVXPY's modeling syntax and uses \texttt{cvxpylayers} to create a differentiable layer. During the forward pass, this layer solves the DC-OPF problem to obtain $\mathbf{p}^{\text{g}}$. In the backward pass, it automatically computes the required Jacobians by differentiating the KKT system, leveraging the methods in~\cite{cvxpylayers2019, amos2017optnet}.}

Next, we will use the collected gradient information defined as $\nabla \mathcal{L}(\mathbf{b}, \boldsymbol{\gamma}, \boldsymbol{\rho})= \mathbf{g}= [{\mathbf{g}^{b}}^\top, {\mathbf{g}^{\gamma}}^\top, {\mathbf{g}^{\rho}}^\top]^{\top}$ in the TNC optimization method to update the coefficient and bias parameters.

\subsection{Optimization Formulation and Solution Method}
\label{sec:Optimization Methods}

To solve \textcolor{black}{the bilevel problem in} Model~\ref{model:O-DCOPF_matrixformat}, we employ the gradient-based TNC optimization method that uses the calculated parameter sensitivities~\cite{nash2000survey, nocedal2006numerical}. \textcolor{black}{It is important to clarify that TNC is used to solve the non-convex upper-level problem~\eqref{eq:optimization problem} of finding the optimal parameters $(\mathbf{b}, \boldsymbol{\gamma}, \boldsymbol{\rho})$. The lower-level DC-OPF problem for each scenario is solved using the Gurobi optimizer, which is invoked through the \texttt{cvxpylayers} library.}
TNC combines the advantages of Newton's method for fast convergence with the conjugate gradient approach to efficiently handle large-scale problems without the need for explicit Hessian computations.
\textcolor{black}{When numerically compared to other common gradient-based optimizers, including Broyden-Fletcher-Goldfarb-Shanno (BFGS) and limited-memory BFGS~\cite{nocedal2006numerical}, we found that TNC consistently demonstrated superior convergence speed and reliability for our specific problem structure, also shown in~\cite{taheri2023optimizing}}.

As shown in Algorithm~\ref{alg:truncated_newton_with_preconditioning}, TNC iteratively updates the solution vector by approximating the Newton direction using truncated conjugate gradients with preconditioning. This begins with an initial guess for the solution vector, followed by preparing an initial diagonal approximation to the Hessian matrix. Each iteration computes the gradient of the loss function and applies the conjugate gradient method to approximate the Newton direction while ensuring that the preconditioner improves numerical stability and accelerates convergence. The stopping criterion for the inner conjugate gradient loop is based on the reduction of the gradient's magnitude. The optimal step length along the search direction is then determined using a line search that satisfies the Wolfe conditions, ensuring sufficient decrease in the objective function. The solution vector is updated accordingly, and the process repeats until the outer convergence criteria are satisfied.

Note that we use a standard implementation of the TNC algorithm to solve the bilevel optimization problem in Model~\ref{model:O-DCOPF_matrixformat}. Other optimization methods could be used without conceptual changes. The novelty of this work lies in the formulation and overall solution approach for the bilevel problem, rather than the application of TNC itself. For a more detailed discussion on the standard TNC method, refer to~\cite{nash2000survey}.

\textcolor{black}{As an alternative solution approach, we could have reformulated the bilevel optimization problem by replacing the lower-level problem with its KKT conditions to obtain a single-level problem~\cite{kleinert2021survey}. Despite its theoretical elegance, this single-level reformulation presents several practical difficulties that can render it intractable. The simultaneous consideration of multiple scenarios yields a large-scale nonlinear program that is computationally demanding. Furthermore, this reformulation introduces non-convexities associated with the bilinear terms from the complementary slackness and stationarity conditions, yielding a mathematical program with equilibrium constraints (MPEC)~\cite{kim2020mpec}. Due to challenges associated with the resulting constraint qualification conditions, this reformulation complicates the solution process, necessitating the use of specialized nonlinear programming solvers whose convergence is notoriously sensitive to various hyperparameters. Thus, for the sake of tractability, we do not explicitly enforce KKT-based constraints. Instead, we follow typical approaches for training machine learning models that rely on sensitivities to obtain gradient-based algorithms that scale well in related settings. Gradients of the loss function with respect to the parameters $(\mathbf{b}, \boldsymbol{\gamma}, \boldsymbol{\rho})$ are obtained by implicitly differentiating the KKT conditions of the lower-level DC-OPF problem, leveraging modern differentiable optimization libraries like \texttt{cvxpylayers}.}

\section{Numerical Experiments}
\label{sec:Numerical experiments}

\begin{figure*}[t!]
\centering
\subfloat[\small $\mathbf{b}$ parameter values]{
    \includegraphics[width=0.95\textwidth]{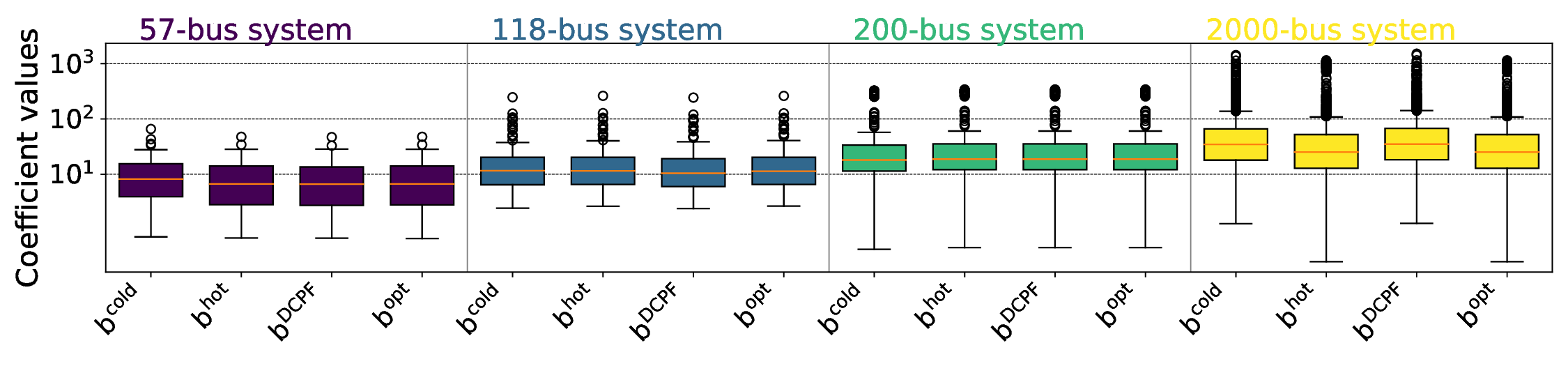}
    \label{fig:box-plot-b}
}
\hfill

\subfloat[\small Injection bias values $\boldsymbol{\gamma}$]{
    \includegraphics[width=0.95\textwidth]{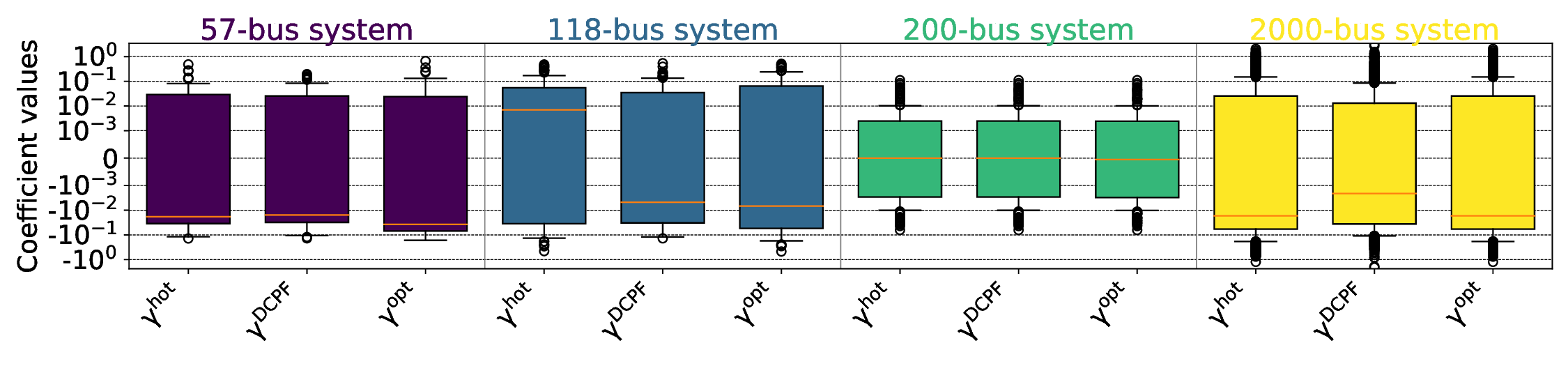}
    \label{fig:boxplot-biases}
}
\hfill

\subfloat[\small Flow bias values $\boldsymbol{\rho}$]{
    \includegraphics[width=0.95\textwidth]{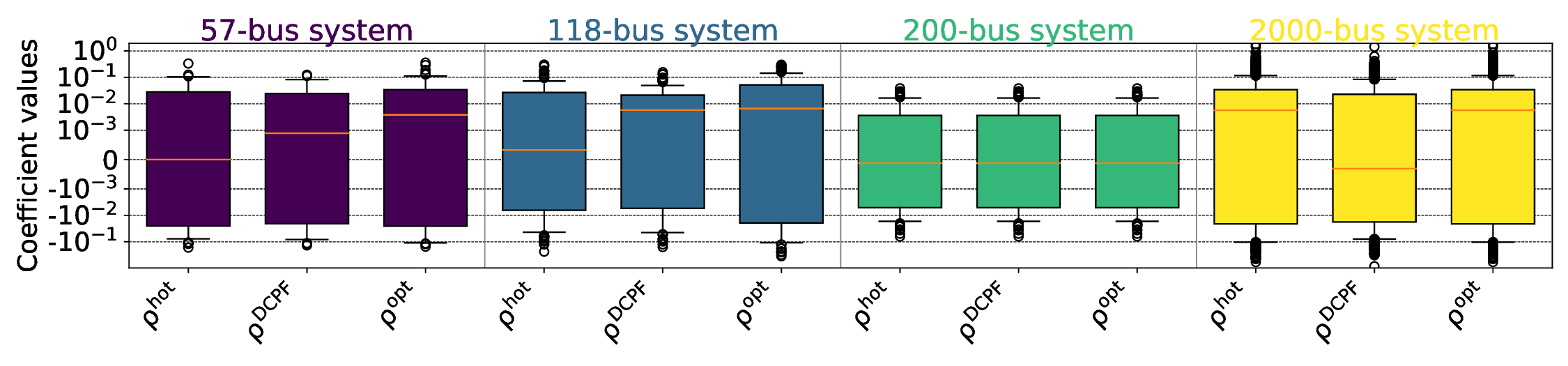}
    \label{fig:boxplot-losses}
}

\caption{Boxplots showing the distributions of the parameter values ($\mathbf{b}$, $\boldsymbol{\gamma}$, and $\boldsymbol{\rho}$) across various test cases for cold-start, hot-start, DCPF, and optimized parameters.}
\label{fig:box-plots}
\end{figure*}

\begin{figure*}[th!]
\centering
\subfloat[\small $\mathbf{b}$ parameter values]{
    \includegraphics[width=0.85\textwidth]{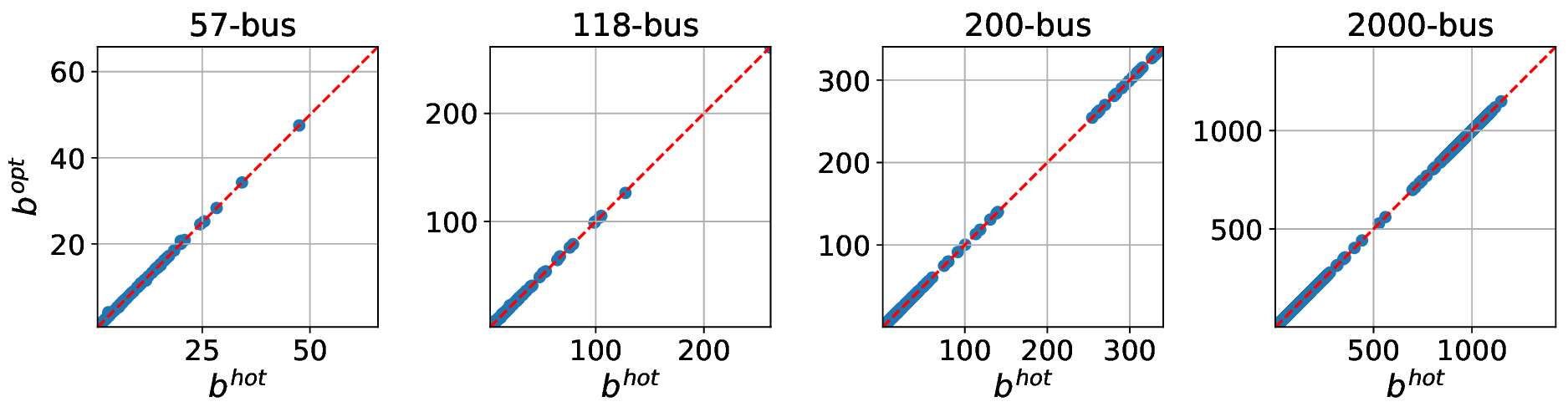}
    \label{fig:scatter-plot-b}
}
\hfill

\subfloat[\small Injection bias values $\boldsymbol{\gamma}$]{
    \includegraphics[width=0.85\textwidth]{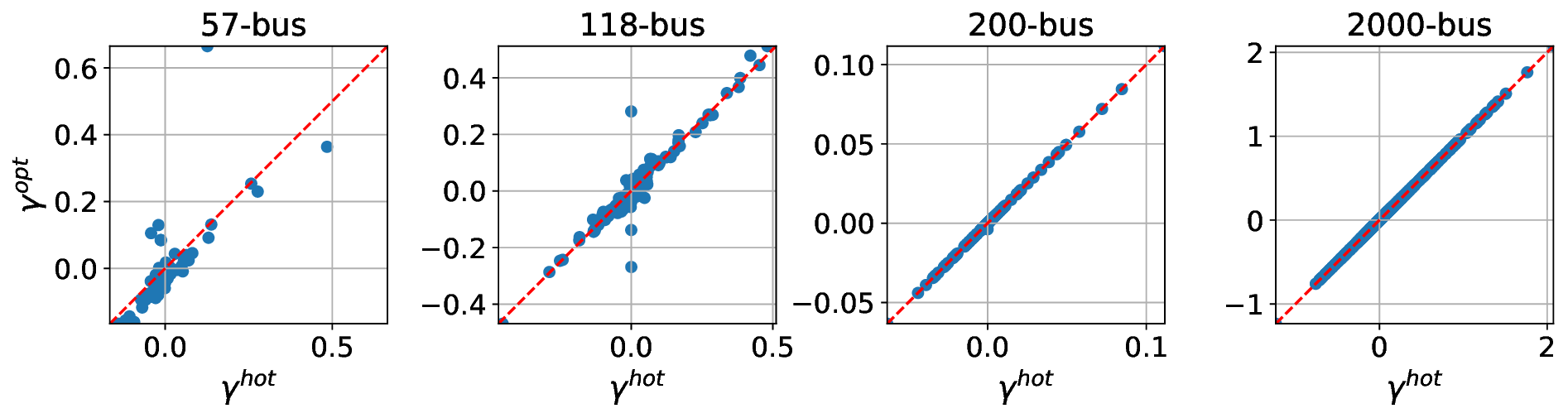}
    \label{fig:scatter-plot-gamma}
}
\hfill

\subfloat[\small Flow bias values $\boldsymbol{\rho}$]{
    \includegraphics[width=0.85\textwidth]{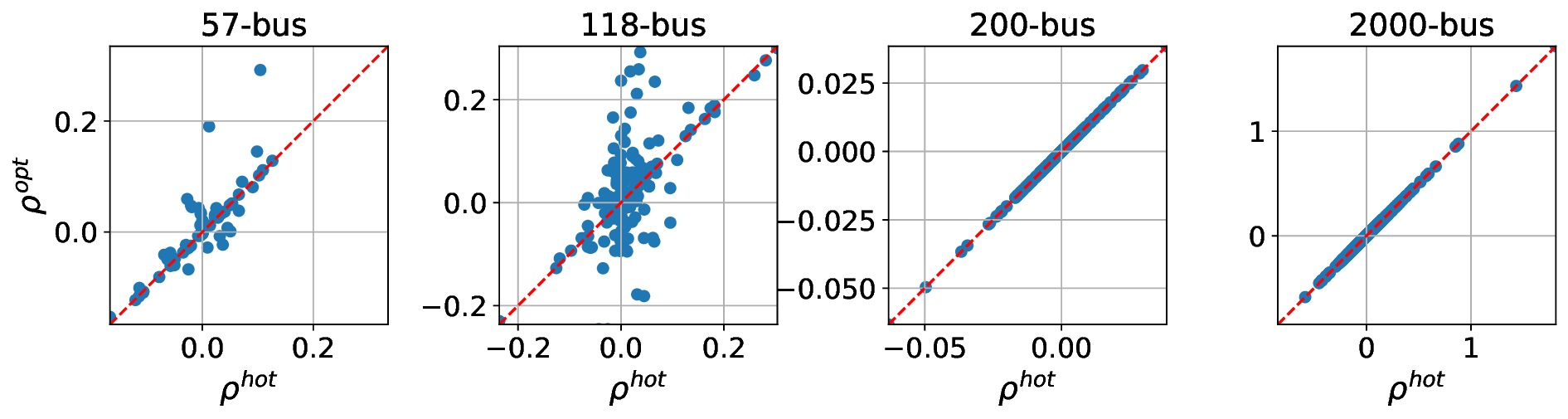}
    \label{fig:scatter-plot-rho}
}
\caption{Scatter plots comparing hot-start ($\mathbf{b}^{hot}$, $\boldsymbol{\gamma}^{hot}$, $\boldsymbol{\rho}^{hot}$) and optimized ($\mathbf{b}^{opt}$, $\boldsymbol{\gamma}^{opt}$, $\boldsymbol{\rho}^{opt}$) parameter values across various test cases. The red dashed line indicates a one-to-one correspondence.}
\label{fig:scatter-plots}
\end{figure*}

Using systems from the PGLib-OPF test case archive~\cite{pglib}, including \textrm{IEEE 14-bus}, \textrm{EPRI 39-bus}, \textrm{IEEE 57-bus}, \textrm{IEEE 118-bus}, \textrm{Illinois 200-bus}, and \textrm{Texas 2000-bus}~\cite{pglib}, this section numerically evaluates our algorithm's performance via comparison against traditional DC-OPF and AC-OPF models.

For each system, we generated $2,020$ load scenarios, allocating $20$ for offline training and $2,000$ for testing, by multiplicatively scaling the nominal loads with normally distributed random variables, each load independently by its own variable, with mean $\mu=1$ and standard deviation $\sigma= 15\%$. The active and reactive demands at each bus were scaled by the same multiplicative factor to keep the loads' nominal power factors.
\textcolor{black}{To ensure a valid basis for comparison, only load scenarios that resulted in a feasible AC-OPF solution were included in the training and testing datasets.
}

All computations were conducted on a 24-core, 32 GB RAM computing node at Georgia Tech's Partnership for an Advanced Computing Environment (PACE). The AC-OPF problems were solved using \textrm{PowerModels.jl}~\cite{coffrin2018powermodels}. We implemented our algorithm in PyTorch with derivatives $\frac{\partial \mathbf{p}^{\text{g}}}{\partial \mathbf{b}}$, $\frac{\partial \mathbf{p}^{\text{g}}}{\partial \boldsymbol{\gamma}}$, and $\frac{\partial \mathbf{p}^{\text{g}}}{\partial \boldsymbol{\rho}}$ calculated via cvxpylayers~\cite{cvxpylayers2019}. The TNC method from the \textrm{scipy.optimize.minimize} library was used to minimize the loss function~\eqref{eq:objective_function} using the sensitivities ${\mathbf{g}^{b}}$, ${\mathbf{g}^{\gamma}}$, and ${\mathbf{g}^{\rho}}$.

We benchmark the accuracy achieved by our optimized parameters by comparing the generators' active power setpoints across different DC-OPF models (i.e., $\mathbf{p}^{[\text{model}]}$ where $[\text{model}]$ stands for either the optimized parameters from our proposed algorithm, the cold-start DC-OPF \eqref{eq:cold-start}, the hot-start DC-OPF~\eqref{eq:hot-start}, or the optimized DC power flow (DCPF) parameters from~\cite{taheri2023optimizing}) against the true values from the \mbox{AC-OPF} solutions (i.e., $\mathbf{p}^{\text{g,AC}}$). \textcolor{black}{For the hot-start DC-OPF benchmark, the parameters are derived from an AC power flow solution at the nominal loads (i.e., the base case with a $1.0$ scaling factor), which represents the center of the distribution used for scenario generation.} The discrepancies are measured using both the mean square error (MSE) as defined in the loss function~\eqref{eq:objective_function} and the maximum error metric, with all values in the per unit (p.u.) system with a $100$~MVA base power:
\begin{align}
\varepsilon_{\text{MSE}}^{[\text{model}]} &= \frac{1}{|\mathcal{G}|\times|\mathcal{M}|} \sum_{m \in \mathcal{M}} \|\mathbf{p}_{m}^{[\text{model}]}-\mathbf{p}_{m}^{\text{g,AC}}\|_{2}^{2},
\end{align}
\begin{align}
    \varepsilon_{\text{max}}^{[\text{model}]} &= \max_{m \in \mathcal{M}} \|\mathbf{p}_{m}^{[\text{model}]}-\mathbf{p}_{m}^{\text{g,AC}}\|_{\infty},
\end{align}
where $\|\cdot\|_{\infty}$ is the $L_{\infty}$-norm and $\|\cdot\|_{2}$ is the $L_{2}$-norm.

\subsection{Comparison of Parameter Values Across Selection Methods}

This section presents the distributions of the cold-start, hot-start, and optimized parameters $\mathbf{b}$, $\boldsymbol{\gamma}$, and $\boldsymbol{\rho}$ for different test systems using box plots and scatter plots. The box plots show the spread of parameter values across the buses and lines in the system, providing insight into how the parameter values vary within the system. The scatter plots compare the hot-start and optimized parameters for each bus or line, demonstrating the consistency of the optimized parameters with those used in a common DC power flow approximation.

\subsubsection{Box Plot Analysis}
Fig.~\ref{fig:box-plots} illustrates the distribution of the parameter values using box plots. The interquartile range (IQR), which represents the middle $50\%$ of the data, is shown in each box, with the central line marking the median value. Whiskers extend to $1.5$ times the IQR, and data points outside this range are considered outliers, marked individually. The box plots for each test system capture the cold-start values ($\mathbf{b}^{cold}$), hot-start values ($\mathbf{b}^{hot}$, $\boldsymbol{\gamma}^{hot}$, and $\boldsymbol{\rho}^{hot}$), and the optimized values ($\mathbf{b}^{opt}$, $\boldsymbol{\gamma}^{opt}$, and $\boldsymbol{\rho}^{opt}$). All values are presented on a logarithmic scale to accommodate variations across multiple orders of magnitude.
The box plots in Fig.~\ref{fig:box-plots} show that the optimized parameter values closely align with those derived from traditional heuristics, particularly for the coefficient parameters. This suggests consistency with longstanding power engineering intuition, where line susceptances are known to be a key factor in determining power flows.
\textcolor{black}{
While the optimization allows negative \(b\) values, all of the optimized $b$ values were positive for all of the test cases. This is consistent with the underlying line reactances in these systems, all of which are positive. If there were lines with negative reactances (e.g., series-compensated lines), the method would yield negative \(b\) values if this would improve~accuracy.
}

\subsubsection{Scatter Plot Analysis}
Complementing the box plots, the scatter plots (Fig.~\ref{fig:scatter-plots}) directly compare the hot-start and optimized parameter values. Each plot includes a $45^\circ$ red dashed line to indicate one-to-one correspondence. Data points falling along this line suggest a close match between the hot-start and optimized values, while deviations highlight discrepancies.
Fig.~\ref{fig:scatter-plot-b} compares the coefficient values $\mathbf{b}^{hot}$ and $\mathbf{b}^{opt}$, while Fig.~\ref{fig:scatter-plot-gamma} compares the injection bias values $\boldsymbol{\gamma}^{hot}$ and $\boldsymbol{\gamma}^{opt}$, and Fig.~\ref{fig:scatter-plot-rho} compares the flow bias values $\boldsymbol{\rho}^{hot}$ and $\boldsymbol{\rho}^{opt}$. These scatter plots reveal that the optimized parameters align well with hot-start values, except for a few outliers where targeted adjustments to $\mathbf{b}$, $\boldsymbol{\gamma}$, and $\boldsymbol{\rho}$ improve the DC-OPF approximation's accuracy.
\textcolor{black}{
While Fig.~\ref{fig:scatter-plots} shows general alignment between optimized and hot-start parameters, we emphasize that this does not imply that the optimized parameters yield marginal accuracy improvements. Hot-start is a strong heuristic, and our algorithm's ability to learn small but systematic deviations from it is precisely what enables substantial accuracy gains across diverse scenarios, as we discuss next.
}

\subsection{Comparison of Approximation Accuracies}

Table~\ref{table:combined_loss} gives the squared two-norm and $\infty$-norm errors for several test cases when using the optimized and conventional DC-OPF parameters over $2,000$ testing scenarios. Notably, in the \textrm{118-bus} test case, our algorithm achieved a substantial reduction in the squared two-norm loss to $0.0123$, compared to the cold-start, hot-start, and optimized DCPF parameters with $0.1295$, $0.1221$, and $0.1251$, respectively. This represents an order of magnitude improvement over all prior methods, highlighting the effectiveness of the proposed algorithm.
\textcolor{black}{
While the absolute \(\infty\)-norm error remains large in some cases, these outliers typically correspond to stressed conditions that any DC power flow approximation struggles to model. Our method cannot fully eliminate all approximation errors, but nevertheless significantly reduces the worst-case error (by as much as $79\%$ across test cases) compared to traditional DC-OPF parameterizations. This shows that the proposed approach substantially increases the DC power flow approximation's accuracy while preserving its computational advantages.}

The improvement column in Table~\ref{table:combined_loss} shows the percentage improvement from the optimized parameters in the DC-OPF models compared to the \textcolor{black}{cold-start, hot-start, and optimized DCPF} parameters. The results show up to $90\%$ and $79\%$ improvements in the squared two-norm and $\infty$-norm errors.

\begin{table*}\centering
\caption{Squared Two-Norm and $\infty$-Norm Loss Functions Using $2,000$ test scenarios}
\label{table:combined_loss}
\setlength{\tabcolsep}{1.5pt} 
\renewcommand{\arraystretch}{1.3}
\resizebox{\textwidth}{!}{%
\begin{tabularx}{\textwidth}{l *{10}{>{\centering\arraybackslash}X} c}
\toprule
\textbf{Test case} & \multicolumn{5}{c}{\textbf{Squared Two-Norm Loss}} & \multicolumn{5}{c}{\textbf{$\infty$-Norm Loss}} \\ 
\cmidrule(lr){2-6} \cmidrule(lr){7-11} 
  & $\varepsilon_{\text{MSE}}^{{cold}}$ & $\varepsilon_{\text{MSE}}^{{hot}}$ & $\varepsilon_{\text{MSE}}^{{DCPF}}$\cite{taheri2023optimizing} &$\varepsilon_{\text{MSE}}^{{opt}}$ &Improv. (\%)& $\varepsilon_{\text{max}}^{{cold}}$ & $\varepsilon_{\text{max}}^{{hot}}$ & $\varepsilon_{\text{max}}^{{DCPF}}$\cite{taheri2023optimizing} &$\varepsilon_{\text{max}}^{{opt}}$&Improv. (\%)& \\ 
\midrule
\midrule
\textrm{14-bus} & 0.0070   & 0.0069   &  0.0069 & $\mathbf{0.0030}$   & (57, 57, 57) & 0.590  & 0.590 & 0.590 & 0.590 & (0, 0, 0)  \\
\textrm{39-bus} & 0.3223  & 0.4248  &  0.3046   & $\mathbf{0.3029}$   & (6, 28, 1) &  5.877  &  6.215  &5.784 & $\mathbf{5.585}$  & (5, 10, 3) \\
\textrm{57-bus} & 0.7445 & 0.6311  & 0.2054& $\mathbf{0.1765}$ & (76, 72, 14) & 4.671   & 4.590 & 3.658&$\mathbf{3.120}$ & (33, 32, 15) \\
\textrm{118-bus}& 0.1295  & 0.1221   & 0.1251& $\mathbf{0.0123}$  & (90, 90, 90) & 3.133   & 3.170 & 3.146 & $\mathbf{1.918}$ & (39, 39, 39) \\

\textrm{200-bus}&0.0031   & 0.0003 & 0.0003 & $\mathbf{0.0002}$  & (93, 33, 33) &0.145  &0.033  & 0.035 & $\mathbf{0.031}$ & (79, 6, 11) \\

\textrm{500-bus}& 0.0073  &  0.0069 & 0.0054 & $\mathbf{0.0030}$  & (59, 56, 44) & 1.356   & 1.000 & 0.950 & $\mathbf{0.890}$ & (34, 11, 6)   \\

\textrm{2000-bus}&  0.0149 & 0.0161  & 0.0156 & $\mathbf{0.0069}$  & (54, 57, 56) &   0.829 &  1.024 & 0.967 & $\mathbf{0.623}$ & (25, 39, 36)   \\

\bottomrule
\end{tabularx}
}
\end{table*}

\begin{figure}[!t]
    \centering
    \includegraphics[width=0.5\textwidth]{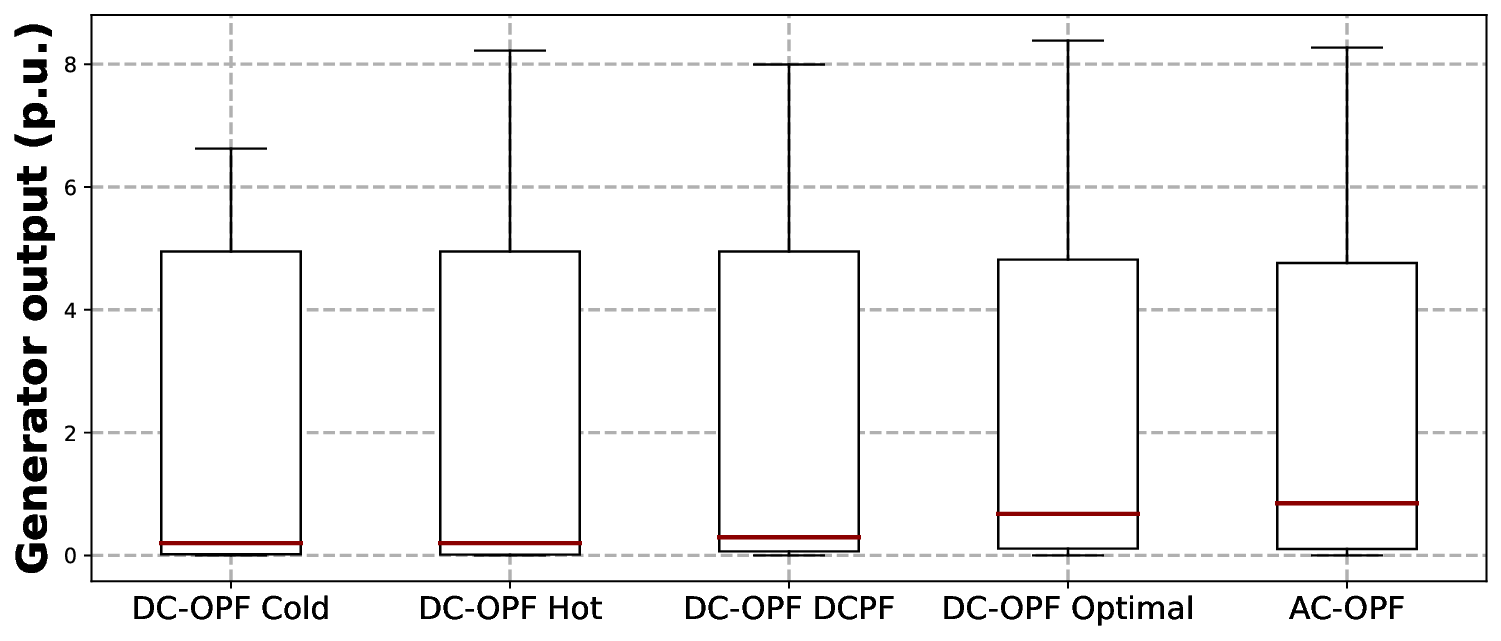}
        \vspace{-1.5em}
    \caption{Comparison of the generator outputs for the \textrm{IEEE 118-bus} system over $2,000$ test scenarios. The box plot compares the performance of the AC-OPF vs DC-OPF with four parameter sets: cold-start ($\mathbf{b}^{cold}$, $\boldsymbol{\gamma}^{cold}$, and $\boldsymbol{\rho}^{cold}$), hot-start ($\mathbf{b}^{hot}$, $\boldsymbol{\gamma}^{hot}$, and $\boldsymbol{\rho}^{hot}$), DCPF-based parameters ($\mathbf{b}^{DCPF}$, $\boldsymbol{\gamma}^{DCPF}$, and $\boldsymbol{\rho}^{DCPF}$), and optimized parameters ($\mathbf{b}^{opt}$, $\boldsymbol{\gamma}^{opt}$, and $\boldsymbol{\rho}^{opt}$).}

    \label{fig:boxplots}
\end{figure}

The box plot in Fig.~\ref{fig:boxplots} compares the mean generator outputs for the \textrm{IEEE 118-bus} system across five formulations: cold-start, hot-start, DCPF, DC-OPF with optimal parameters, and AC-OPF. The vertical axis represents the generator output in per unit, and the horizontal axis shows the different parameter choices. Each box illustrates the distribution of mean generator outputs for non-zero generators across $2,000$ test scenarios, highlighting the variability and median for each method. The figure shows that the DC-OPF with optimized parameters ($\mathbf{b}^{opt}$, $\boldsymbol{\gamma}^{opt}$, and $\boldsymbol{\rho}^{opt}$) aligns closely with the AC-OPF, as indicated by similar medians and ranges.

\begin{figure*}[!t]
\centering
\subfloat[\small Generator setpoints]{
    \includegraphics[width=0.48\textwidth]{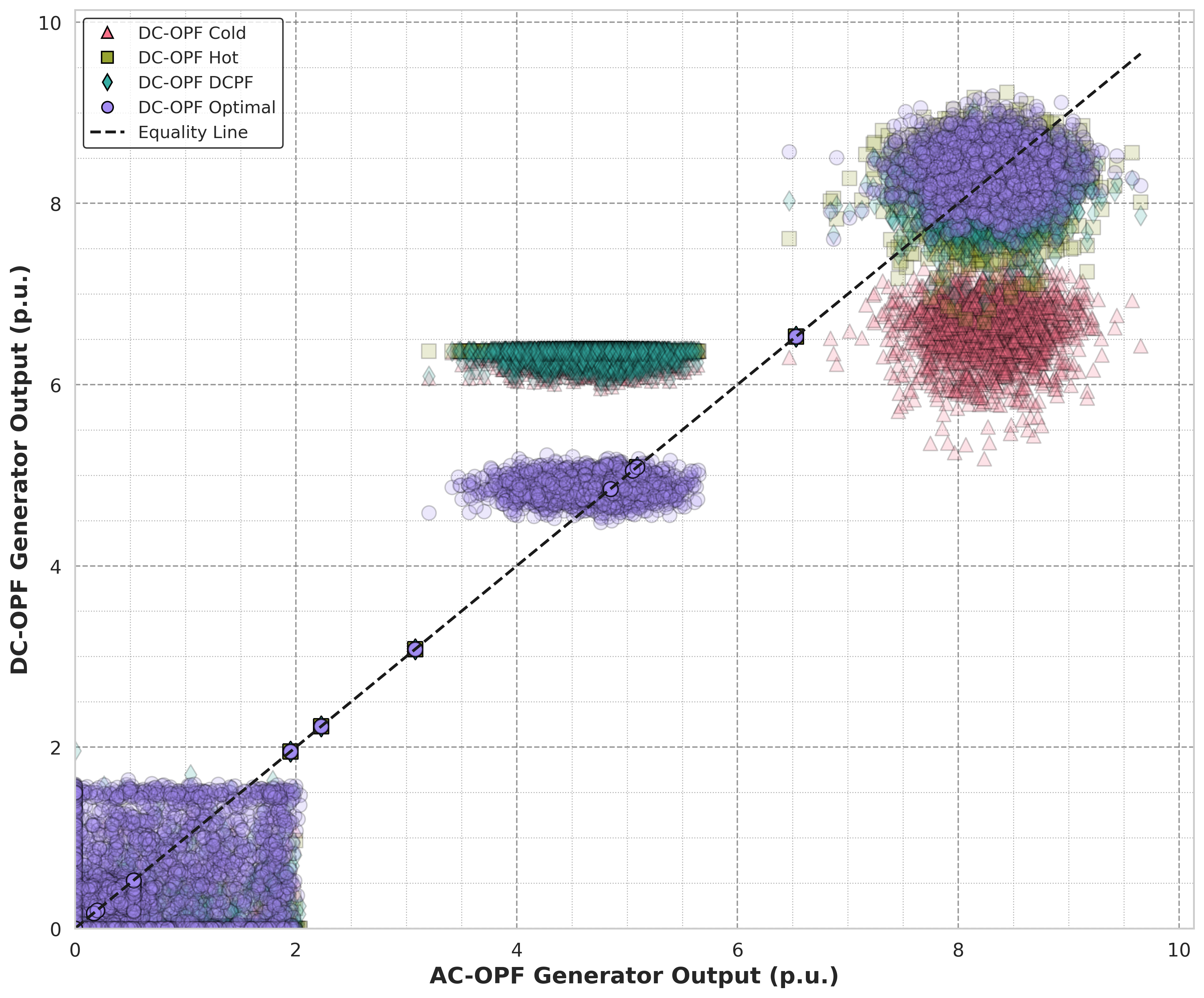}
    \label{fig:generator outputs1}
}
\hfill
\subfloat[\small Average generator setpoints]{
    \includegraphics[width=0.48\textwidth]{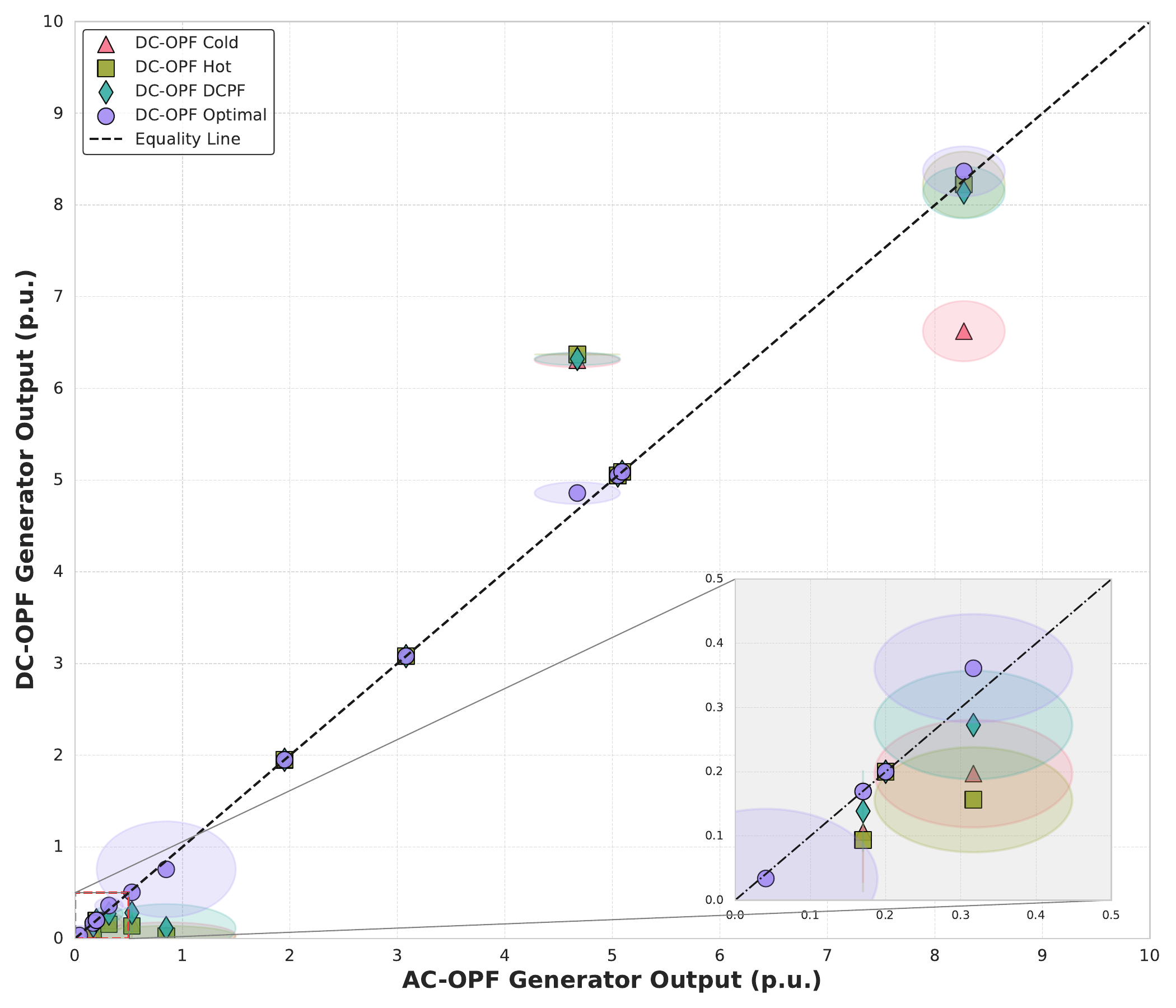}
    \label{fig:generator outputs2}
}
\vspace{-0.5em}
\caption{(a) Generator setpoints (p.u.) obtained by solving the AC-OPF problem and DC-OPF problems with cold-start, hot-start, DCPF, and optimized parameters for the \textrm{IEEE 118-bus} test system. Generators that consistently reached their maximum limits across all models are excluded. The 45-degree equality line indicates perfect alignment between the AC-OPF and DC-OPF setpoints, serving as a benchmark for comparison. (b) Average generator setpoints (p.u.) with standard deviation ellipses across 2,000 scenarios.}
\label{fig:generator outputs}
\end{figure*}

\begin{figure*}[!th]
    \centering
    \includegraphics[width=1\textwidth]{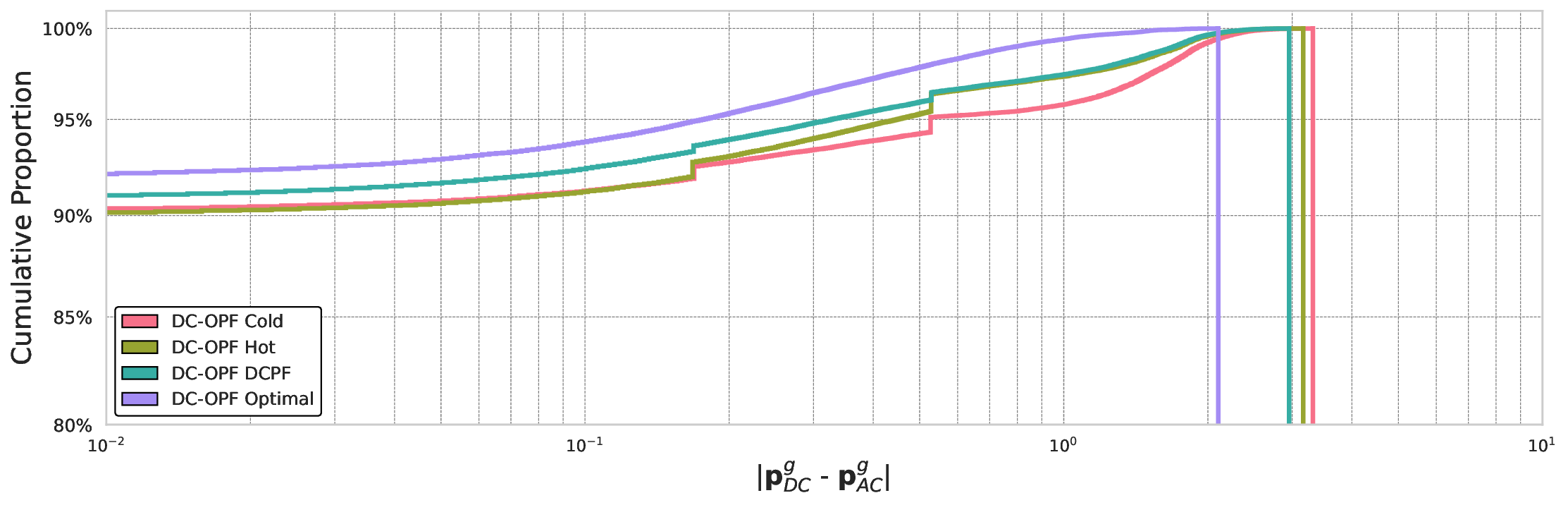}
    \vspace{-2em}
    \caption{Cumulative proportion of the absolute error between AC-OPF and DC-OPF (in per unit) for the \textrm{IEEE 118-bus} system over $2,000$ test scenarios. The graph compares four scenarios: usage of cold-start ($\mathbf{b}^{cold}$, $\boldsymbol{\gamma}^{cold}$, and $\boldsymbol{\rho}^{cold}$), hot-start ($\mathbf{b}^{hot}$, $\boldsymbol{\gamma}^{hot}$, and $\boldsymbol{\rho}^{hot}$), DCPF ($\mathbf{b}^{DCPF}$, $\boldsymbol{\gamma}^{DCPF}$, and $\boldsymbol{\rho}^{DCPF}$), and optimized ($\mathbf{b}^{opt}$, $\boldsymbol{\gamma}^{opt}$, and $\boldsymbol{\rho}^{opt}$) parameters, shown on a logarithmic scale. }
    \label{fig:cumulative}
\end{figure*}

The graphical analyses provided in Figs.~\ref{fig:generator outputs} and~\ref{fig:cumulative} offer further insights into the performance of the proposed algorithm. In Fig.~\ref{fig:generator outputs}, the generator setpoints for the \textrm{IEEE 118-bus} test system are presented, where we compare the AC-OPF solutions (i.e., $\mathbf{p}^{\text{g,AC}}$) with the DC-OPF solutions derived from models with cold-start, hot-start, DCPF, and optimized DC-OPF parameters. Fig.~\ref{fig:generator outputs1} displays the setpoints of individual generators, excluding those that consistently reached their maximum output limits across all parameter choices. The alignment of setpoints with the 45-degree equality line demonstrates how closely the DC-OPF models, particularly the one with optimized parameters, mimic the AC-OPF results. 
Fig.~\ref{fig:generator outputs2} shows the average generator setpoints (p.u.) with standard deviation ellipses across $2,000$ scenarios. These ellipses illustrate the variability and performance of each DC-OPF method, highlighting the consistency of the optimized DC-OPF in closely matching the AC-OPF setpoints. The inset focuses on generators with lower setpoints, providing a detailed view of their alignment and variability. These results demonstrate that our optimized DC-OPF parameters achieve setpoints comparable to the more complex AC-OPF while offering enhanced computational efficiency.

Fig.~\ref{fig:cumulative} shows the cumulative proportion of the absolute errors between AC-OPF and DC-OPF solutions for the \textrm{IEEE 118-bus} system. The graph compares four scenarios: usage of cold-start $\mathbf{b}$, hot-start $\mathbf{b}^{hot}$, $\boldsymbol{\gamma}^{hot}$, $\boldsymbol{\rho}^{hot}$, DCPF $\mathbf{b}^{DCPF}$, $\boldsymbol{\gamma}^{DCPF}$, $\boldsymbol{\rho}^{DCPF}$, and optimized $\mathbf{b}^{opt}$, $\boldsymbol{\gamma}^{opt}$, $\boldsymbol{\rho}^{opt}$ parameters, shown on a logarithmic scale. Points located towards the upper left represent better performance, indicating a larger proportion of smaller errors, which demonstrates the improved accuracy of the DC-OPF solutions using the optimized parameters.

Observe that the maximum errors are less than $1.918$ per unit for the solution obtained from the optimized DC-OPF model versus errors up to $3.133$, $3.170$, and $3.146$ per unit resulting from the cold-start DC-OPF with $\mathbf{b}^{cold}$, the hot-start DC-OPF with $\mathbf{b}^{hot}$, $\boldsymbol{\gamma}^{hot}$, and $\boldsymbol{\rho}^{hot}$, and the DCPF with $\mathbf{b}^{DCPF}$, $\boldsymbol{\gamma}^{DCPF}$, and $\boldsymbol{\rho}^{DCPF}$, respectively.

\begin{table}[t]
    \centering
    \caption{CPU Time in Seconds}
        \label{table:CPU time}
    \renewcommand{\arraystretch}{1.3} 
    \label{tab:comparison}
    \resizebox{\columnwidth}{!}{\Large
    \begin{tabular}{lccccccc}
        \toprule
        
        \textbf{Test case} & \textrm{14-bus} & \textrm{39-bus}   & \textrm{57-bus}  & \textrm{118-bus} &  \textrm{200-bus} & \textrm{500-bus} & \textrm{2000-bus} \\
        \midrule
        \midrule
        $t_{\text{train}}$ & 13 & 17 & 41 & 52  &58 & 903 & 9357 \\
        \hdashline
        $t_{\text{DCOPF}}$ & 0.009 & 0.021 & 0.024 & 0.064 & 0.055& 0.354 & 1.304 \\
        $t_{\text{ACOPF}}$ & 0.030 & 0.071 & 0.091 & 0.295 & 0.289& 1.650 & 8.708 \\
        \bottomrule
    \end{tabular}
    }
    \begin{tablenotes}
        \item[*] \scriptsize $t_{train}$: Offline computation time to train the parameters.
        \item[*] \scriptsize $t_{DCOPF}$: Online per-scenario computation time averaged over $100$ scenarios.
         \item[*] \scriptsize $t_{ACOPF}$: Online per-scenario computation time averaged over $100$ scenarios.
    \end{tablenotes}
\end{table}


Finally, Table~\ref{table:CPU time} summarizes the training, testing, and AC-OPF model computation times. Training times, denoted as $t_{train}$, range from $13$ to $9,357$ seconds. This represents a one-time, upfront effort that is within reasonable ranges for offline computations (minutes to hours for larger systems). We note that the primary bottlenecks in training computations are solving the DC-OPF problems and computing the gradients of the loss function with respect to the parameters. On average, solving the DC-OPF problems accounts for approximately $45\%$ of the total training time, while computing the gradients contributes to $55\%$ of the training time. We note that an implementation that utilizes a stochastic gradient descent method, as is often employed when training machine learning models, could reduce both components of the training time.

Post-training, the solution times for the DC-OPF model ($t_{DCOPF}$) range from $0.009$ to $1.304$ seconds per scenario, significantly faster than the AC-OPF model due to the simplicity of solving linear programs versus nonlinear programs. Moreover, the variation in DC power flow parameters does not impact DC-OPF solution times, which remain consistent to within 2\% across different parameter choices.

\textcolor{black}{While the times in Table~\ref{table:CPU time} show that a single AC-OPF solve is computationally acceptable in some contexts, the key advantages of a fast and accurate DC-OPF model emerge in more complex applications. Problems such as unit commitment, optimal transmission switching, or stochastic optimization often involve solving many thousands of OPF subproblems~\cite{padhy2004, barrows2014, roald2022review}. For these applications, AC-OPF computational costs are often prohibitive, making a high-accuracy DC-OPF model indispensable. The proposed algorithm provides parameters that enhance DC-OPF accuracy for these computationally intensive tasks without sacrificing its essential speed advantage.}

\textcolor{black}{
Compared to our prior work optimizing parameters to match DC and AC \textit{power flows}~\cite{taheri2023optimizing}, the bilevel approach in this paper is computationally more intensive since differentiating through a constrained optimization problem (DC-OPF) is more complex than differentiating through the closed-form power flow equations. The benefits of the method's enhanced accuracy outweigh its additional computational cost. As shown in Table~\ref{table:combined_loss}, directly optimizing for the end-goal—accurate generator setpoints from an OPF problem—yields significantly greater accuracy improvements (e.g., a $90\%$ vs. $14\%$ improvement in MSE for the \textrm{57-bus} case when comparing the proposed method to the DCPF method of~\cite{taheri2023optimizing}). This superior performance justifies the added one-time offline computational effort for applications where solution quality is paramount.}

\section{Conclusion}
\label{sec:conclusion}
This paper presented an optimization algorithm for significantly enhancing the accuracy of DC-OPF problems by selecting coefficient and bias parameters in the DC power flow approximation based on a sensitivity analysis. Employing the TNC method for optimization alongside efficient gradient computation techniques, an offline training phase refines DC power flow parameters across a range of loading scenarios. With squared two-norm and \(\infty\)-norm losses that are up to $90\%$ and $79\%$, respectively, lower than those of traditional parameter selection methods for DC-OPF problems, our numerical tests demonstrate a substantial improvement in DC-OPF solution accuracy relative to the AC-OPF problem.

Future work will evaluate the optimized parameters in complex applications like unit commitment and infrastructure hardening. We will also test their generalization to security-constrained problems, including $N-1$ contingencies, and their robustness against more diverse daily and seasonal load profiles.
\textcolor{black}{
Our preliminary experience suggests that achieving the most substantial improvements will likely require adapting the bilevel formulation and loss function to the specific objectives and constraints of each application (e.g., discrete decisions and multi-period constraints in unit commitment).}
\textcolor{black}{
As the complexity of the training data increases, for example, by including generator commitment variations, the optimization landscape may become more challenging. While TNC proved effective in this study, future implementations can easily adopt other optimizers through the same \texttt{scipy.optimize} interface to assess their performance in other applications. Options such as stochastic gradient-based or hybrid methods are standard for large-scale, non-convex problems and could further improve the speed and reliability of convergence under more realistic and complex operational scenarios.}
Finally, we plan to analyze the underlying system characteristics, such as network congestion and topology, to determine the conditions under which our method provides the most significant accuracy gains.

\bibliographystyle{IEEEtran}
\bibliography{blib}

\begin{IEEEbiography}[{\includegraphics[width=1in, height=1.25in, clip, keepaspectratio]{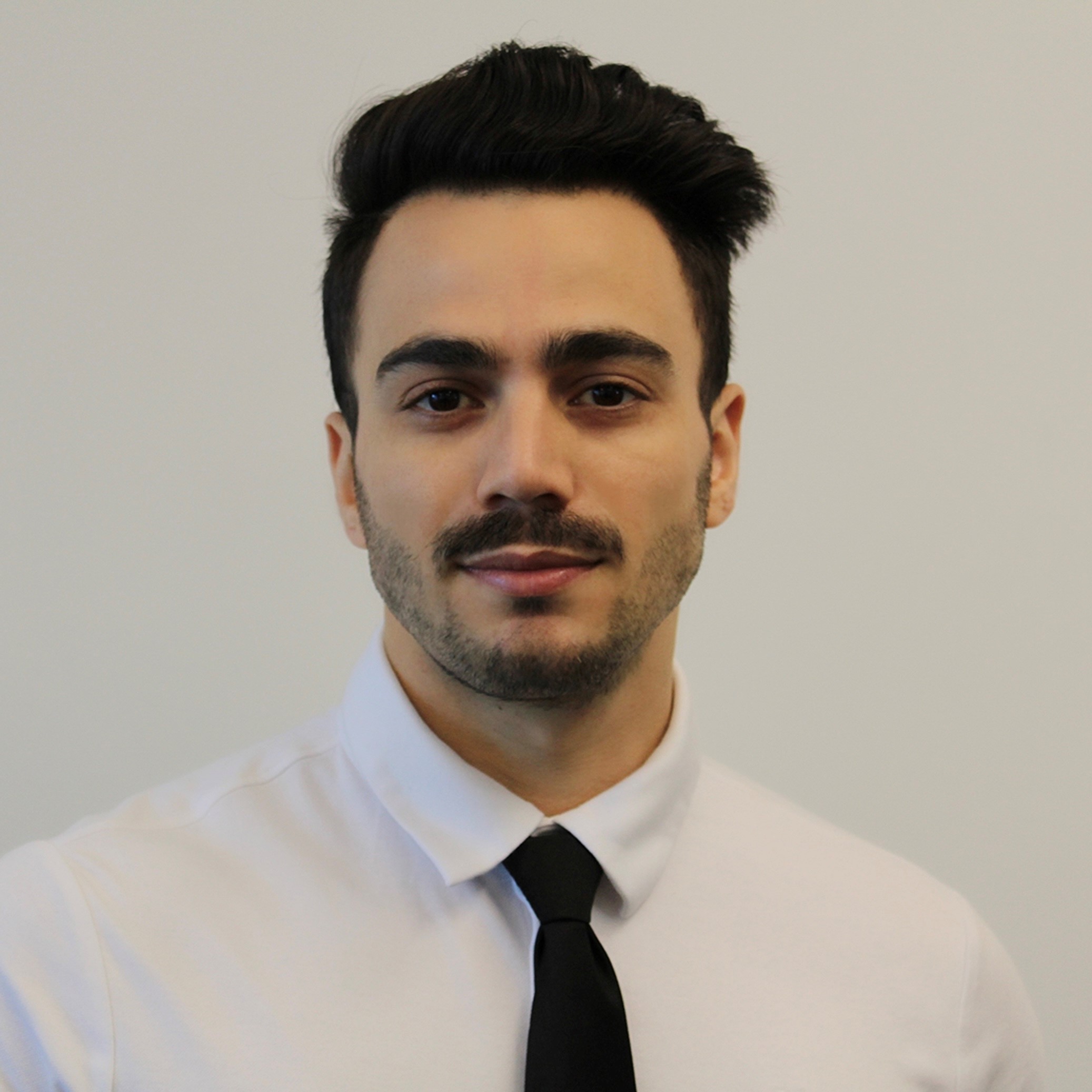}}]{Babak Taheri} received the Ph.D. and M.Sc. degrees in electrical and computer engineering from the Georgia Institute of Technology, Atlanta, GA, USA in 2024. He also received the M.Sc. degree in electrical engineering from Sharif University of Technology in 2019 and the B.Sc. degree from the University of Tabriz in 2017. He is currently a Research Scientist at Hitachi Energy Research. His research interests include power systems, optimization, and machine learning.
\end{IEEEbiography}

\vspace{-9pt}

\begin{IEEEbiography}[{\includegraphics[width=1in, height=1.25in, clip, keepaspectratio]{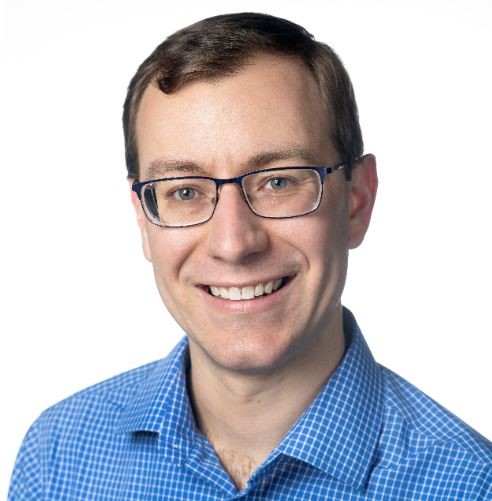}}]{Daniel K. Molzahn} (Senior Member, IEEE) received
the B.S., M.S., and Ph.D. degrees in electrical engineering and the master’s of Public Affairs degree from the University of Wisconsin--Madison, Madison, WI, USA. He is currently an Associate Professor with the School of Electrical and Computer Engineering, Georgia Institute of Technology, Atlanta, GA, USA, and also holds an appointment as a Computational
Engineer in the Energy Systems Division at Argonne National Laboratory. He was a Dow Postdoctoral Fellow in Sustainability at the University of Michigan, Ann Arbor, MI, USA, and a National Science Foundation Graduate Research Fellow at the University of Wisconsin--Madison. He was the recipient of the IEEE Power and Energy Society’s Outstanding Young Engineer Award in 2021, the NSF CAREER Award in 2022, and Georgia Tech’s Class of 1940 W. Roane Beard Outstanding Teacher Award in 2024.
\end{IEEEbiography}

\end{document}